\documentclass{article}

\usepackage[]{arxiv}
\PassOptionsToPackage{hyphens}{url}
\usepackage[numbers,sort&compress]{natbib}%
\bibpunct[, ]{[}{]}{,}{n}{,}{,}%
\usepackage[table,xcdraw]{xcolor} %
\usepackage{scalerel}
\usepackage{tikz,setspace}
\usetikzlibrary{svg.path}
\definecolor{orcidlogocol}{HTML}{A6CE39}
\tikzset{
  orcidlogo/.pic={
    \fill[orcidlogocol] svg{M256,128c0,70.7-57.3,128-128,128C57.3,256,0,198.7,0,128C0,57.3,57.3,0,128,0C198.7,0,256,57.3,256,128z};
    \fill[white] svg{M86.3,186.2H70.9V79.1h15.4v48.4V186.2z}
                 svg{M108.9,79.1h41.6c39.6,0,57,28.3,57,53.6c0,27.5-21.5,53.6-56.8,53.6h-41.8V79.1z M124.3,172.4h24.5c34.9,0,42.9-26.5,42.9-39.7c0-21.5-13.7-39.7-43.7-39.7h-23.7V172.4z}
                 svg{M88.7,56.8c0,5.5-4.5,10.1-10.1,10.1c-5.6,0-10.1-4.6-10.1-10.1c0-5.6,4.5-10.1,10.1-10.1C84.2,46.7,88.7,51.3,88.7,56.8z};
  }
}
\newcommand\orcidicon[1]{\href{https://orcid.org/#1}{\mbox{\scalerel*{
\begin{tikzpicture}[yscale=-1,transform shape]
\pic{orcidlogo};
\end{tikzpicture}
}{|}}}}
\usepackage{graphicx}
\usepackage[utf8]{inputenc}
\usepackage{amsfonts}
\usepackage{tabularx,booktabs}

\usepackage{amsmath}
\usepackage{mathtools}
\usepackage{bbm}
\usepackage{makecell}
\newcommand{\argmin}{\operatornamewithlimits{argmin}}

\usepackage{amsmath,amssymb,amsfonts}
\usepackage{textcomp,setspace}
\usepackage{float}
\usepackage[utf8]{inputenc} %
\usepackage[T1]{fontenc}    %
\usepackage{nth}
\usepackage{algorithm,algpseudocode}       %
\usepackage{url}            %
\usepackage{booktabs}       %
\usepackage{amsfonts}       %
\usepackage{nicefrac}       %
\usepackage{microtype}      %
\usepackage{array}
\usepackage{amsbsy} 
\usepackage{amsmath}
\usepackage{ulem}
\usepackage{enumitem}
\usepackage{blkarray}
\usepackage{amsmath}

\def\BibTeX{{\rm B\kern-.05em{\sc i\kern-.025em b}\kern-.08em
		T\kern-.1667em\lower.7ex\hbox{E}\kern-.125emX}}
\makeatletter
\def\BState{\State\hskip-\ALG@thistlm}
\makeatother

\newcommand{\stt}{{\rm s.t.}}
\newcommand{\cG}{\mathcal{G}}
\newcommand{\cV}{\mathcal{V}}
\newcommand{\cE}{\mathcal{E}}

\newcommand{\cN}{\mathcal{N}}

\newcommand{\cR}{\mathcal{R}}

\newcommand{\br}{\mathbf{r}}
\newcommand{\bR}{\mathbf{R}}
\newcommand{\bb}{\mathbf{b}}
\newcommand{\bB}{\mathbf{B}}
\newcommand{\bT}{\mathbf{T}}
\newcommand{\bv}{\mathbf{v}}
\newcommand{\bs}{\mathbf{s}}
\newcommand{\bS}{\mathbf{S}}

\newcommand{\bc}{\mathbf{c}}
\newcommand{\bq}{\mathbf{q}}
\newcommand{\bD}{\mathbf{D}}

\newcommand{\bu}{\mathbf{u}}
\newcommand{\bH}{\mathbf{H}}
\newcommand{\bW}{\mathbf{W}}
\newcommand{\bI}{\mathbf{I}}

\newcommand{\bZ}{\mathbf{Z}}
\newcommand{\bgamma}{\boldsymbol{\gamma}}
\newcommand{\bbeta}{\boldsymbol{\beta}}

\newcommand{\bla}{\boldsymbol{\lambda}}
\newcommand{\bLa}{\boldsymbol{\Lambda}}
\newcommand{\bomega}{\boldsymbol{\omega}}
\newcommand{\balpha}{\boldsymbol{\alpha}}
\newcommand{\bdelta}{{\boldsymbol{\delta}}}
\newcommand{\bDelta}{{\boldsymbol{\Delta}}}
\newcommand{\bmu}{\boldsymbol{\mu}}
\newcommand{\bbR}{\mathbb{R}}

\usepackage{eso-pic}
\usepackage{multirow}
\usepackage{times}
\usepackage{fancyhdr,graphicx,amsmath,amssymb}
\include{pythonlisting}
\usepackage{etoolbox, lineno}
\usepackage{lipsum}

\usepackage{soul}
\usepackage[]{color-edits}
\addauthor{mr}{blue}
\addauthor{ali}{red}
\usepackage[breaklinks]{hyperref}%
\begin{document}
\title{Incentive Systems for Fleets of New Mobility Services}

\author{Ali~Ghafelebashi \orcidicon{0000-0001-8339-7960}, Meisam~Razaviyayn \orcidicon{0000-0003-4342-6661}, and~Maged~Dessouky \orcidicon{0000-0002-9630-6201} \thanks{This study was funded by a grant from the National Center for Sustainable Transportation (NCST), supported by the U.S. Department of Transportation’s University Transportation Centers Program. The contents of this project reflect the views of the authors, who are responsible for the facts and the accuracy of the information presented herein. This document is disseminated in the interest of information exchange. The U.S. Government and the State of California assume no liability for the contents or use thereof. Nor does the content necessarily reflect the official views or policies of the U.S. Government and the State of California. This paper does not constitute a standard, specification, or regulation. This paper does not constitute an endorsement by the California Department of Transportation (Caltrans) of any product described herein.\newline Ali Ghafelebashi, Meisam Razaviyayn, and Maged Dessouky are with the Daniel J. Epstein Department of Industrial \& Systems Engineering, University of Southern California, 3715 McClintock Avenue, Los~Angeles, CA 90089, United States.\newline \textit{Email addresses:} ghafeleb@usc.edu (Ali Ghafelebashi), razaviya@usc.edu (Meisam Razaviyayn), maged@usc.edu (Maged Dessouky).}}

\maketitle

\begin{abstract}
Traffic congestion has become an inevitable challenge in large cities due to population increases and expansion of urban areas.
Various approaches are introduced to mitigate traffic issues, encompassing from expanding the road infrastructure to employing demand management.
Congestion pricing and incentive schemes are extensively studied for traffic control in traditional networks where each driver is a network ``player''.
In this setup, drivers' ``selfish'' behavior hinders the network from reaching a socially optimal state.
In future mobility services, on the other hand, a large portion of drivers/vehicles may be controlled by a small number of companies/organizations.  
In such a system,  offering incentives to organizations can potentially be much more effective in reducing traffic congestion rather than offering incentives directly to drivers. 
This paper studies the problem of offering incentives to organizations to change the behavior of their individual drivers (or individuals relying on the organization’s services). 
We developed a model where incentives are offered to each organization based on the aggregated travel time loss across all drivers in that organization. 
Such an incentive offering mechanism requires solving a large-scale optimization problem to minimize the system-level travel time. 
We propose an efficient algorithm for solving this optimization problem. 
Numerous experiments on Los~Angeles County traffic data reveal the ability of our method to reduce system-level travel time by up to 6.9\%.
Moreover, our experiments demonstrate that incentivizing organizations can be up to 8 times more efficient than incentivizing individual drivers in terms of incentivization monetary cost.
\end{abstract}

\keywords{New Mobility Services \and Congestion Reduction \and Incentivizing Organizations \and Travel Demand Management}

\section{Introduction}
\label{sec:introduction}
Today, traffic congestion is one of the major issues in metropolitan areas across the globe. Traffic congestion declines the overall quality of life, leads to significant economic losses, degrades air quality, and escalates health vulnerabilities due to emissions~\cite{INRIX, air_polution, health2010traffic, zhang2013air}. This paper undertakes the task of devising a novel mechanism around incentives. The core objective of this mechanism is to change the behavioral patterns of individual drivers within organizations by incentivizing organizations. 

\vspace{0.2cm}

Incentive-based congestion reduction methodologies overlap with pricing methods, including taxes and fees for road access~\cite{pigou1920economics, knight1924some, RePEc:elg:eebook:4192, van2009behavioural, bouchelaghem2018reliable, zhang2013self, kachroo2016optimal, farokhi2014study, groot2014toward, cao2020improving}. These strategies encourage individuals to avoid congested routes, reducing traffic buildup. Various determinants underpin the design of these pricing frameworks, encompassing temporal aspects~\cite{zheng2016time}, spatial metrics~\cite{daganzo2015distance}, and vehicular attributes~\cite{zhang2019impact, zhong2020pricing}. Although promising, market-oriented pricing and taxation face challenges due to equity concerns, policy complexity, and implementation uncertainties~\cite{knockaert2012spitsmijden, levinson2010equity, martens2012justice, raux2004acceptability, ieromonachou2006evaluation, hensher2014type, gu2018congestion}.

\vspace{0.2cm}

Another approach within the area of pricing mechanisms involves the adoption of tradable credits (TCs) or tradable mobility permits (TPMs)~\cite{verhoef1997tradeable, wang2014models, fan2013tradable, thogersen2008breaking}. \cite{tsekeris2009design} provides a theoretical analysis of the benefits of tradable credits. This methodology has been implemented within some economic sectors, exemplified by its use in the airport slot market~\cite{fukui2010empirical}. Nevertheless, the implementation of these cap-and-trade programs in personal travel and daily commutes is hindered by design challenges~\cite{dogterom2017tradable, azevedo2018tripod}.

\vspace{0.2cm}

Recently, there has been a heightened focus on incentivization strategies. Compared to fee-based methods, reward-based policies can be more popular~\cite{brehm1966theory}. 
Moreover, the efficacy of incentivizing positive actions over punishing negative ones is evidenced in the psychological concept of reactance~\cite{brehm1966theory}.
While rewarding policies have proven effective in altering individual behavior~\cite{kreps1997intrinsic, berridge2001reward}, the transportation sector has underexplored these incentives.

\vspace{0.2cm}

There have been several studies that explored the use of incentivization to reduce traffic congestion, such as the INSTANT project~\cite{merugu2009incentive}, the CAPRI project of Stanford~\cite{yue2015reducing}, series of studies in the Netherlands~\cite{bliemer2009rewarding}, and the ``Metropia'' platform~\cite{hu2015behavior}. 
In a recent study, \cite{wu2023managing} showed the effectiveness of ridesharing incentivization in congestion reduction. 
Incentivizing off-peak hour driving is examined via public and private central planners (policymakers) in \cite{wu2023public}. 
However, congestion reduction by offering incentives to organizations has not been studied by any of the previous studies. 
Although initial success is shown in reward-based strategies, enduring behavior change is not always guaranteed~\cite{kumar2016impacts}.

\vspace{0.2cm}

\begin{figure*}[] 
  \centering
  \begin{tabular}{@{}c@{}}
    \includegraphics[width=1.\linewidth]{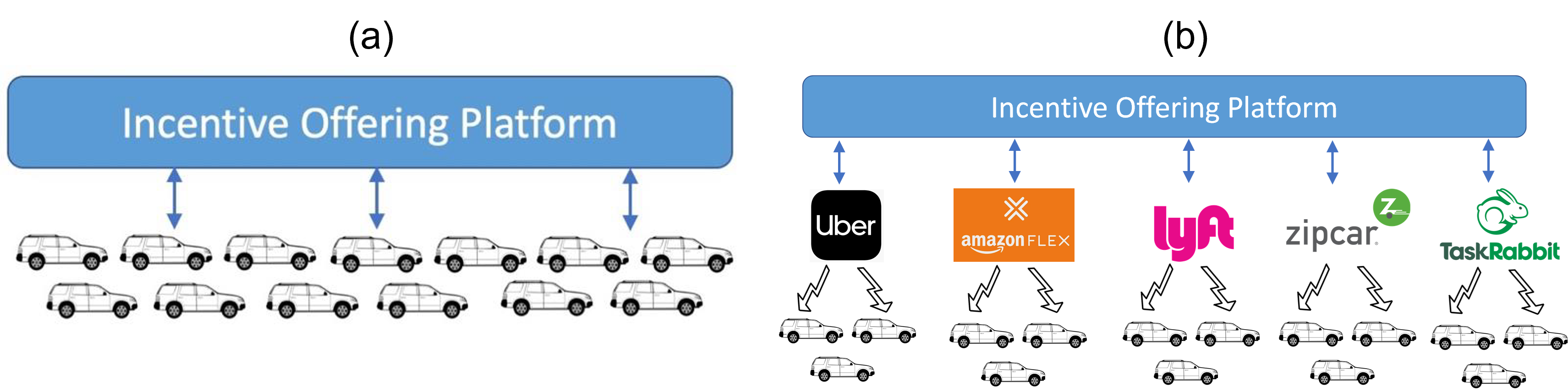} \\[\abovecaptionskip]
  \end{tabular}  
  \vspace{-0.50cm}
  \caption{(a) Traditional platforms for offering incentives. (b) Presented platform for offering incentives.}
  \label{fig:incentiveOfferingPlatforms}
\end{figure*}

In traditional congestion pricing and incentive offering mechanisms, incentives are offered directly to individual drivers to influence their decisions, such as departure time and routing (Figure~\ref{fig:incentiveOfferingPlatforms} (a)). In mobility services, many of these decisions may be directly (or indirectly) made by organizations providing different transportation services. For example, navigation apps, which are regularly used by almost 70\% of smartphone users~\cite{he2019people, Panko}, influence the routing decision of millions of drivers daily.
Another example is crowdsourcing delivery platforms such as Amazon Flex, Instacart, and Doordash. According to a recent study~\cite{DoorDashGrowth}, the revenue of DoorDash during the fourth quarter of 2022 increased by 40\% to \$1.8 billion from \$1.3 billion in revenue that it recorded during the same period in 2021.
Another example of such organizations is ride-hailing companies such as Uber and Lyft. According to a report by Uber for the fourth quarter of 2022~\cite{UberFinance2022}, the number of gross bookings increased from 12\% in the fourth quarter of 2021 to 17.7\% in the fourth quarter of 2022. Today, many of the routing decisions are made by individual drivers. With the future emergence of autonomous vehicles, it is possible that organizations may now own the fleet of vehicles and control their routing. 
Intuitively, since organizations have more flexibility and more power to change the traffic, incentivizing organizations is expected to be more efficient than incentivizing individual drivers.  Furthermore, an organization has more options in balancing the route selection across the large pool of drivers employed by the organization. Motivated by this idea,  this paper develops an incentive offering mechanism to organizations to indirectly (or directly) influence the behavior of individual drivers  (Figure~\ref{fig:incentiveOfferingPlatforms} (b)).
In a different study, \cite{ghafelebashi2023congestion} utilizes the traditional incentive offering framework (Figure~\ref{fig:incentiveOfferingPlatforms} (a)) to provide an algorithm to offer personalized incentives to drivers to reduce traffic congestion by changing the routing decision of the drivers. These incentives could be personalized based on user preferences.

\vspace{0.2cm}

In contrast to the individual-level incentivization presented in \cite{merugu2009incentive, yue2015reducing, bliemer2009rewarding, ghafelebashi2023congestion}, our incentivization framework for organizations addresses a broader spectrum of challenges and complexities:
\begin{enumerate}
    \item  
    Numerical experiments in~\cite{ghafelebashi2023congestion} show that  the monetary benefit from reduced travel time based on the Value of Time (VOT) exceeds the incentivization cost. 
    However, their model does not depend on VOT. Our model uses VOT to compute the monetary value of organizations' time loss and compensate for it through incentivization.
    
    \item Some incentivization studies~\cite{merugu2009incentive, yue2015reducing, bliemer2009rewarding} offered static rewards based on fixed rules for all the participants. However, our model utilizes different VOTs to compute the incentive offer for different organizations. Note that different organizations can have different VOTs due to the unique nature of their service. Hence, our model offers incentives to organizations such that they can compensate the organizations' time loss based on their VOT.
    
    \item \cite{ghafelebashi2023congestion} offers personalized incentives, but they are selected from a discrete set of incentive choices that are fixed before solving the problem. However, our model employs a continuous variable to calculate the value of the required incentive. This variable depends on VOT and the amount of drivers' time loss. As we are not limited to a discrete set of incentives, our incentivization cost can be more cost-efficient. Note that the variability in incentive values introduces more complexity to our optimization problem because of the larger variable size. 
    
    \item \cite{ghafelebashi2023congestion} does not consider fairness and time delivery constraints (fair assignment of drivers to slower and faster routes when they share the same origin and destination simultaneously). 
    In contrast, our model addresses these limitations by preventing the diversion of drivers to routes with significant time disparities, thus ensuring a route assignment based on fairness and time delivery constraints.  
\end{enumerate}

\vspace{0.2cm}

Our framework will be based on the following three-step procedure:

\noindent Step~1) The central planner receives organizations’ demand estimates for the next time interval (e.g., the next few hours).

\noindent Step~2) The central planner incentivizes organizations to change their routes and travel time.

\noindent Step~3) Observe organizations’ response and go back to Step~1 for the next time interval.

\noindent The central planner (which is referred to as “Incentive Offering Platform” in Figure~\ref{fig:incentiveOfferingPlatforms} (b)) continually repeats this three-step process in the network for every time interval.

\vspace{0.2cm}

The rest of this paper is structured as follows. Section~\ref{sec:whyOrganizations} motivates the advantage of incentivizing organizations compared to incentivizing individuals. Section~\ref{sec:IncentiveOfferingMechanism} introduces the basic notations and describes our incentive offering mechanism for congestion reduction. We formulate an optimization problem to find the ``optimal'' incentive offering strategy. We then propose an algorithm for solving this optimization problem efficiently in Section~\ref{sec:ADMMAlgo}. Numerical experiment results for the model using Los Angeles County data are detailed in Section \ref{sec:NumericalExperiments}. Concluding remarks are discussed in Section~\ref{sec:conclusion}

\vspace{-0.2cm}

\section{Why Offering Incentives to Organizations Rather Than Individuals?}  \label{sec:whyOrganizations}
\begin{figure}[H] 
  \centering
  \begin{tabular}{@{}c@{}}
    \includegraphics[width=0.45\linewidth]{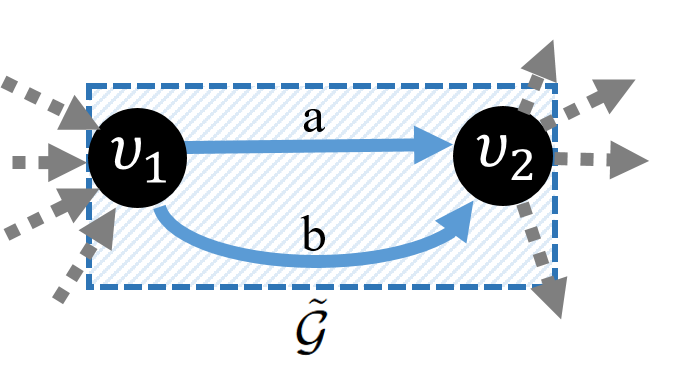} \\[\abovecaptionskip]
  \end{tabular}  
  \vspace{-0.50cm}
  \caption{Subnetwork $\tilde{\mathcal{G}}$ (selected in blue dashed rectangle).}
  \label{fig:toyExampleNetwork}
\end{figure}

\vspace{-0.2cm}
Our methodology is incentivizing organizations (rather than individual drivers). Let us first motivate the benefit of this strategy via a simple example. Consider the subnetwork $\tilde{\mathcal{G}}$ at Figure~\ref{fig:toyExampleNetwork} as a subset of a larger network.
\vspace{-0.55cm}
\vspace{-0.3cm}
Links $a$ and $b$ are routes between nodes $v_1$ and $v_2$. The travel times of $a$ and $b$ are 25 and 30 minutes at User Equilibrium (UE), respectively. Assume 20 drivers start traveling from $v_1$ to $v_2$ at the same time. If travel time is the most important factor in their utility, they will select $v_1$ because it is the fastest route at UE. Assume we have found the System Optimal (SO) strategy for the entire network, and we need 15 out of the 20 drivers to select $b$ instead of $a$ to achieve SO. At SO, the travel time of route $a$ decreases to 20 minutes (5-minute decrease), and the travel time of route $b$ increases to 35 minutes (5-minute increase). Drivers that use route $a$ save 5 minutes due to a decrease in travel time of route $a$. Deviated drivers to route $b$ expect to lose 5 minutes because they expect route $b$ to have travel time of 30 minutes. Hence,  since we want to deviate 15 drivers to a route with longer travel time (route $b$ in this example), we should compensate for their increased travel time. Assume VOT is $\$1/\text{min}$. Let us  compare two scenarios:
\begin{enumerate}[label=(\Roman*)]
    \item All 20 drivers are individual drivers. Since we need to deviate/incentivize 15 drivers and compensate each of them for $5$ minutes of their time, we need to spend $\$75 = (5\; \text{min} \times 15)\times \$1/\text{min}$. \label{ScenraioI} 
     \item All $20$ drivers work in the same organization. In this scenario, the organization needs to spend \$75 to alter the decision of the $15$ drivers. However, after offering incentives, the travel time of the 5 remaining drivers on route~$a$ decreases. Therefore, the organization gains $25 = 5 \times 5$ minutes of time from the drivers who stayed in route~$a$. Overall, the increase and decrease in the travel times of the drivers cancel each other out (canceling-out effect). This change only costs the organization 50~minutes of total time. Hence, the compensation cost is $\$50 = 50\; \text{min} \times \$1/\text{min}$ for the organization. \label{ScenraioII}
\end{enumerate}
Therefore, we spend 33\% less in incentivizing the organization (i.e., scenario \ref{ScenraioII}) compared to incentivizing the individual drivers (i.e., scenario \ref{ScenraioI}). The above simple example illustrates that incentivizing organizations can be more cost-effective than incentivizing individual drivers. 
Note that this observation does not necessarily hold in general games. That is, grouping users in a game does not necessarily lead to a lower-cost Nash equilibrium.

\vspace{0.2cm}

\section{Incentive Offering Mechanism and Problem Formulation } \label{sec:IncentiveOfferingMechanism}
\vspace{-0.2cm}
Given the origin-destination information of drivers in various organizations, the goal is to find the ``optimal'' strategy for offering organization-level incentives to them to reduce the traffic congestion of the system. To mathematically state the problem, we begin this section by defining our notations. A complete list of notations used in this paper can be found in Appendix~\ref{apdx:notations}. For further details of the notation, an example is provided in Appendix~\ref{appdx:ModelNotationExample}. 

\vspace{0.2cm}

The traffic network is represented by a directed graph~$\cG = (\cV,\cE)$. Vertices~$\cV$ of the graph are major ramps and intersections in the network. Vertices are connected by a set of edges~$\cE$. 
In our directed graph, the edge direction is determined by the allowable direction of travel on the corresponding road segment, indicating the permissible movement from one node to another for a driver.
The adjacency of two nodes is based on the possibility of driving directly from one node to another without visiting any other node. 
The network comprises a total number of road segments, denoted as~$|\cE|$, which reflects the cardinality of the set~$\cE$.
A route in the network is a path in the graph and is denoted by a one-hot encoding. In other words, a given route is represented by a vector $\mathbf{r} \in \{0,1\}^{|\cE|}$ in which the $k$-th entry is one if route~$\mathbf{r}$ includes the $k$-th edge and it is zero, otherwise. Let $\bT = \{1,\ldots,T\}$ denote the defined time horizon such that $t=1$ marks the starting time of the system. Traffic volume of road segments at time~$t$ is represented by the vector $\bv_t \in \mathbb{R}^{|\cE|}$ in which the $k$-th entry is the total number of vehicles of road segment~$k$ at time~$t$.

\vspace{0.2cm}

Let $\cN = \cN_1\cup\dots\cup\cN_n$ denote the set of all drivers and $\cN_i $ denote the set of drivers of organization~$i$. If a driver works for multiple organizations, he or she will be counted as a different driver at each organization. Hence, $\cN_1 \cap \dots \cap \cN_n = \emptyset$. For any driver $j \in \cN$, let $\cR_j \subseteq \{0,1\}^{|\cE|}$ denote the set of driver's possible route choices between her/his origin and destination. The binary variable $s_i^{\br,j}\in\{0,1\}$ represents the assigned route to the~$j$-th driver of organization~$i$. For this driver and given route~$\br \in \cR$, the variable $s_i^{\br,j} = 1$ if route~$\br$ is assigned to the~$j$-th driver of organization~$i$; and $s_i^{\br,j} = 0$, otherwise. Each driver can only be assigned to one route, i.e., $\sum_{\br \in \cR_j} s_i^{\br,j} = 1$. Given any routing strategy assigned to drivers, we model the drivers' decision deterministically due to the power of the organizations in enforcing their drivers' routes. 

\vspace{0.2cm}

In this paper, we change the routing decision of organizations' drivers by incentivizing their organizations. 
We assume that organizations will accept our route assignments if the incentive offer can compensate for the change in their total travel time. Notice that when the organizations decide to accept the offer, they have no access to the offered route assignments to the other organizations. Hence, they can only estimate the travel time based on historical data, and they will be compensated based on their loss/gain compared to the historical setting.

\vspace{0.2cm}

In this work, we adopt total travel time as the utility function, while alternative metrics like energy consumption or total carbon emissions can also be considered.
We compute the system total travel time by summing the drivers' travel time of all road segments over all time periods in the horizon of interest:
\begin{equation} \label{ObjFunct1}
    \begin{split}
        F(\hat{\bv}) = & \sum_{\ell=1}^{|\cE|} \sum_{t=1}^{|\mathbf{T}|} \hat{v}_{\ell, t} \theta_{\ell, t}(\hat{v}_{\ell, t})  \\
    \end{split}
\normalsize
\end{equation}
\normalsize
where $\theta_{\ell, t}$ is the travel time of link~$\ell$ at time~$t$ (which itself is a function of the link's traffic volume at that specific time). Here, $\hat{\bv}$ is the vector of volume of links in which $\hat{v}_{\ell, t}$ is the $(|\cE|\times (t-1) + \ell)^{th}$ element of vector~$\hat{\bv}$ corresponding to the volume on the~$\ell^{th}$ link at time~$t$. Using the volume vector, we can then calculate the travel time for the links at various times, as outlined below.

\vspace{0.2cm}

Multiple approaches have been proposed to illustrate the relationship between traffic volume and travel time. For instance, the Bureau of Public Roads (BPR) \cite{united1964traffic} presents a congestion function for road links. This function describes a non-linear connection between the travel time on a road and its traffic volume:
\begin{equation} \label{BPR}
    \begin{aligned}
       \theta(v) =  f_{\text{BPR}}(v) = \theta_0 \left(1+0.15\left(\frac{v}{w}\right)^4\right)
    \end{aligned}
\normalsize
\end{equation}
\normalsize
where $f_{\text{BPR}}(v)$ denotes the travel time for drivers on a road segment based on its traffic volume~$v$; $\theta_0$ represents the segment's free flow travel time; and $w$ is the road segment's practical capacity.
In our experiments, to learn the parameters $w$ and $\theta_0$ of the road segments in the Los~Angeles area at different times of the day, we utilize the historical traffic data of the road segments. Given the function $\theta(\cdot)$ in~\eqref{BPR}, to compute the total travel time of the system, one needs to compute the volume at each link. Subsequently, we elucidate the process by which the volume vector is computed within our model.

\vspace{0.2cm}

\noindent\textit{Volume vector~$\hat{\bv}$:} The computation of the volume vector~$\hat{\bv}$ requires (approximately) estimating the location of the drivers at different times based on their route. By assigning a different route to a driver, the driver's impact on the values of the vector~$\hat{\bv}$ will be different because the driver's location will change by following a different route. We will begin by introducing our notation for route assignment: Each driver's assigned route is represented by a one-hot encoded vector. Thus, for each driver, we have a binary vector $\bs_i^j \in \{0, 1\}^{|\cR|}$ in which only one element has a value of one, and it corresponds to the assigned route to the $j$-th driver of organization~$i$.  As we need one vector for each driver, we can aggregate all our assignments in a matrix $\bS \in \{0, 1\}^{|\cR| \times |\cN|} = \left[\bS_1 \bS_2 \dots \bS_n\right]$ where $\bS_i \in \{0, 1\}^{|\cR| \times |\cN_i|}$, which is the assignment matrix of organization~$i$ with $n$ being the number of organizations. 
Elements in a driver's assignment vector that correspond to routes unrelated to their specific origin-destination pair are set to zero since travel on these routes is not possible for the drivers.
Thus, the $\bS$ matrix can be rather sparse.

\vspace{0.2cm}

Given the driver's route entering the system at a specific time, we need to model the location of the individual in the upcoming times. To model the drivers' location in the system, we use the model developed by \cite{ma2018estimating} in which the location of drivers is computed in a probabilistic fashion. This model can be presented by a matrix~$\bR  \in [0, 1]^{(|\cE|\cdot|\bT|) \times |\cR|} $ which estimates the probability of a driver being on a certain link at a given future time, under the assumption that they choose a specific route.  Multiple ways to estimate matrix~$\bR$ are suggested in \cite{ma2018estimating}, including an approach based on the use of historical data. In our experiments in subsection~\ref{sec:simModel}, matrix~$\bR$  is computed based on the volume at the UE state of the system. 
Given matrix~$\bR$, it is easy to see that the vector
\begin{equation} \label{volumeVector}
    \begin{split}
        \hat{\bv} =  & \bR \bS \mathbf{1} \in \bbR^{|\cE|\cdot|\bT|} \\
    \end{split}
\normalsize
\end{equation}
\normalsize
contains the expected number of vehicles in all the links at each time. Plugging the expression of $\hat{\bv}$ in \eqref{ObjFunct1}, we get the total travel time of the system as
\begin{equation} \label{eq:F_tt}
    \begin{split}
        F(\hat{\bv}) & = \sum_{\ell=1}^{|\cE|} \sum_{t=1}^{|\mathbf{T}|} (\bR\bS\mathbf{1})_{\ell, t} \theta((\bR\bS\mathbf{1})_{\ell, t}) \\
        & = \sum_{\ell=1}^{|\cE|} \sum_{t=1}^{|\mathbf{T}|} (\br_{\ell, t}\bS\mathbf{1}) \theta(\br_{\ell, t}\bS\mathbf{1}) 
    \end{split}
    \normalsize
\end{equation}
\normalsize
where $\br_{\ell, t}$ is the row of matrix $\bR$ which corresponds to link $\ell$  at time $t$. 

\vspace{0.2cm}

To reduce the total travel time of the system, some drivers can be deviated to alternative routes to lower the traffic flow of the congested links. To change the routing assignment of drivers, we need to offer incentives to their organizations such that it can compensate the organizations' financial loss caused by accepting our assignment. For simplicity, we use the total travel time increase of the organization as a measure of financial loss. Although we have estimated the travel time of the system from equation~\eqref{eq:F_tt}, we need to compute the ``route travel times'' to be able to compare the amount of change in travel time of each driver after offering incentives. 
Given the route travel times, we compute the incentives using a model that depends on VOT and the amount of increase in the travel time for each organization. In particular, we assume that, given the route assignment to organization~$i$, the incentive value is 
\begin{equation} \label{eq:IncentiveValue}
    \begin{split}
c_{i} = \alpha_i \max\left\{0,\sum_{j \in\cN_i} \bdelta^{\top} \bs_i^{j} - \gamma_i\right\},
    \end{split}
\normalsize
\end{equation}
\normalsize
where $c_{i}$ is the incentive offered to organization~$i$, $\alpha_i   \in \bbR_{+}$ is VOT for a driver in organization~$i$, $\bdelta \in \bbR_{+}^{|\cR|\cdot|\bT|} $ is the travel time of the route for each driver, 
and $\gamma_i$ is the sum of the minimum travel time route of each driver of organization~$i$ in the absence of incentivization. $\bdelta$ and $\gamma_i$ are computed based on the absence of incentivization. When $\sum_{j \in\cN_i} \bdelta^{\top} \bs_i^{j} - \gamma_i>0$, the organization's total travel time has increased compared to the baseline of having no incentive, and hence the system will compensate the organization's loss. On the other hand, when  $\sum_{j \in\cN_i} \bdelta^{\top} \bs_i^{j} - \gamma_i<0$, the organization's travel time is improved after incentivization, and hence no incentivization is required for this particular organization to participate.  The details of our method for computing route travel time vector~$\bdelta$ are described next. 

\vspace{0.2cm}

\noindent\textit{Route travel time vector~$\bdelta$:} Estimation of the vector~$\bdelta$ requires the volume on each link which is derived based on the route assignment of the drivers. Let $\bS$ denote the routing decision of the drivers. Given $\bS$, we can estimate the volume vector~$\bv$ using (\ref{volumeVector}). By utilizing the BPR function~(\ref{BPR}) and the estimated volume vector~$\bv$, we can compute the speed of the links. Given the speed of each link, we can determine the vector~$\bdelta$ that contains the travel time of the routes for different time slots and the vector~$\boldsymbol{\eta} \in \bbR_{+}^{K\cdot|\bT|}$ that contains the travel time of the fastest route for different OD pairs for different times. ($K$ represents the total number of origin-destination (OD) pairs). To do so, we rely on the method provided by \cite{ma2018estimating} and the routing decision of drivers $\bS$ at the UE state of the system.
Given the minimum travel time between OD pairs in vector~$\boldsymbol{\eta}$, we can compute the minimum travel time of organization~$i$ as $\gamma_i = (\bB_i \boldsymbol{\eta})^{\top} \mathbf{1}$ where  $\bB_i \in \{0, 1\}^{|\cN_i| \times (K \cdot |\bT|)}$ is the matrix of shortest travel time assignment of drivers of organization~$i$. $\bB_i \boldsymbol{\eta}$ is the vector of the shortest travel time between the OD pair for each driver, and by summing the elements of this vector, we get $\gamma_i$. 

\vspace{0.2cm}

\normalsize\noindent\textit{Proposed formulation:} For minimizing the total travel time of the system via providing incentives to organizations,
we need to solve the following optimization problem: 
\begin{equation} \label{eq:optModel1}
    \begin{split}
        \min_{\{\bS_i, c_i\}_{i=1}^{n}} \quad & \sum_{\ell=1}^{|\cE|} \sum_{t=1}^{|\mathbf{T}|} \hat{v}_{\ell, t} \theta_{\ell, t}(\hat{v}_{\ell, t}) \\ 
        \stt \quad & \hat{\bv} = \sum_{i=1}^{n} \bR \bS_i \mathbf{1} \\
        & \bD \bS_i \mathbf{1} = \bq_i,\quad \forall i=1, 2, \dots, n \\
        & \bS_i^{\top} \mathbf{1} = \mathbf{1},\quad  \forall i=1, 2, \dots, n \\
        & \bS_i \in \lbrace0, 1\rbrace^{(|\cR|\cdot|\bT|)\times(|\cN_{i}|)},\quad  \forall i=1, 2, \dots, n \\
        & \bS_i^{\top} \bdelta \leq \bb_i \odot \bB_i \boldsymbol{\eta},\quad  \forall i=1, 2, \dots, n\\
        & c_i \geq \alpha_i(\bdelta^{\top} \bS_i \mathbf{1} - \gamma_i),\quad \forall i=1, 2, \dots, n \\
        & c_1 + c_2 + \dots + c_n \leq \Omega \\
        & c_i \geq 0,\quad \forall i=1, 2, \dots, n \\ 
    \end{split}
\end{equation}        
\normalsize
where $\hat{v}_{\ell, t}$ is an element of vector~$\hat{\bv}$ that corresponds to the volume of link $\ell$ at time $t$, $c_i \in \bbR_{+}$ is the cost of incentive assigned to organization~$i$, $\bD \in \{0, 1\}^{(K \cdot |\bT|) \times (|\cR| \cdot |\bT|)}$ is the matrix of route assignment of the OD pairs, $\bb_i \in \bbR_{+}^{|\cN_i|}$ denotes the factor by which the travel time of an assigned route can be larger than shortest travel time of the OD pair, $\bB_i \in \{0, 1\}^{|\cN_i| \times (K \cdot |\bT|)}$ is the matrix of shortest travel time assignment of drivers of organization~$i$, and $\bq_i \in {\bbR}^{K \cdot |\bT|}$ is the vector of the number of drivers of organization~$i$ for each OD pair at different times. If there are drivers in the system that do not work for any organization, we can consider them as a single organization whose decision matrix is initialized and has fixed values such that they are assigned to the fastest route (assuming they always select the shortest route). The same idea can be employed for organizations not joining the incentivization platform. The following section provides a detailed explanation of the constraints:

\vspace{0.2cm}

\noindent\textbf{Constraint 1} ($\hat{\bv} = \sum_{i=1}^{n} \bR \bS_i \mathbf{1}$): This constraint is the computation of the volume on each link at different times based on the routing assignments for the organizations. 

\vspace{0.2cm}

\noindent\textbf{Constraint 2} ($\bD \bS_i \mathbf{1} = \bq_i$): This constraint ensures that the correct number of drivers are assigned to the routes between OD pairs. $\bS_i \mathbf{1}$ represents the number of drivers that have been assigned to the different routes. 
The matrix $\bD$ is utilized to aggregate the count of drivers assigned to various routes within the same origin-destination (OD) pair. 
The vector $\bq_i$ represents the actual number of drivers from organization~$i$ traveling between these OD pairs, and the product $\bD \bS_i \mathbf{1}$ is required to equal $\bq_i$.

\vspace{0.2cm}

\noindent\textbf{Constraint 3} ($\bS_i^{\top} \mathbf{1} = \mathbf{1}$): This constraint simply states that we can only assign one route to each driver of organization~$i$. 

\vspace{0.2cm}

\noindent\textbf{Constraint 4} ($\bS_i \in \{0, 1\}^{(|\cR|\cdot|\bT|) \times |\cN_i|}$):  
This constraint enforces a binary framework on our decision variables, where $0$ indicates not assigning a route and $1$ signifies route assignment.

\vspace{0.2cm}

\noindent\textbf{Constraint 5} ($\bS_i^{\top} \bdelta \leq \bb_i \odot \bB_i \boldsymbol{\eta}$): This is our fairness and time delivery constraint. Due to different reasons, such as urgent deliveries by some of the organizations' drivers, they may not accept alternative routes that deviate significantly from the fastest route. Moreover, the platform should consider fairness between different drivers in terms of the amount of deviation from the shortest travel time. The fairness and time delivery constraint bounds the deviation of travel time of the assigned routes from the minimum travel time. $\bS_i^{\top} \bdelta$ represents the travel time of the assigned routes to drivers of organization~$i$.  $\bb_i \in \bbR_{+}^{|\cN_i|}$ denotes the factor by which deviation is allowed for each driver.

\vspace{0.2cm}

\noindent\textbf{Constraint 6} ($c_i \geq \alpha_i(\bdelta^{\top} \bS_i \mathbf{1} - \gamma_i)$ and $c_i \geq 0$):  These two constraints guarantee~\eqref{eq:IncentiveValue}.

\vspace{0.2cm}

\noindent\textbf{Constraint 7} ($c_1 + c_2 + \dots + c_n \leq \Omega$): This represents our budget constraint. 
The scalar $c_i$ denotes the incentive amount allocated to organization~$i$. $\Omega$ signifies the total budget available.

\vspace{0.2cm}

For further elaboration on model~\ref{eq:optModel1} and its constraints, an illustrative example is presented in Appendix~\ref{appdx:ModelNotationExample}.

\vspace{0.2cm}

\section{Incentivization Algorithm and A Distributed Implementation} \label{sec:ADMMAlgo}
\vspace{-0.1cm}
Optimization problem~\eqref{eq:optModel1} is of large size and includes binary variables ($\bS_i, \forall i=1,\dots,n$). Thus, solving it efficiently is a challenging task. 
In this subsection, we propose an efficient algorithm for solving it. First, we  relax the binary constraint~$\bS_i \in \{0, 1\}^{(|\cR|\cdot|\bT|) \times |\cN|}$ to convex constraint~$ \bS_i \in [0, 1]^{(|\cR|\cdot|\bT|) \times |\cN|}$ and we refer to this as the relaxed version of problem~\eqref{eq:optModel1}.    
The objective function is a summation of monomial functions with positive coefficients. Furthermore, $\theta_{\ell, t}$ is an affine mapping of the optimization variable $\bS_i$. Since our domain is the nonnegative orthant and monomials are convex in this domain, the objective function is convex. As the constraints of this problem are convex, the relaxed version of problem~\eqref{eq:optModel1} becomes a convex optimization problem. 
Thus, standard solvers such as CVX~\cite{cvx} can be used to solve this problem. However, these solvers have large computational complexity because of utilizing methods such as interior point methods~\cite{gb08} with $O(n^3)$ iteration complexity where $n$ is the number of variables. This computational complexity is not practical for our problem. In what follows, we rely on first-order methods with linear computational complexity in $n$, which is affordable in our problem. 
The reformulation is provided in Appendix~\ref{appndx:optModelReformulatedSection}. As we discuss in Appendix~\ref{appdx:ADMM}, this reformulation is amenable to the ADMM method~\cite{boyd2011distributed,gabay1976dual, glowinski1975approximation,hong2016convergence, barazandeh2021efficient}, which is a first-order method and scalable. The steps of the resulting algorithm are provided in Algorithm~\ref{alg:ADMM-ForOurProblem} in Appendix~\ref{appdx:alg}. The details of the derivation of this algorithm are provided in Appendix~\ref{appdx:ADMMSteps}. Due to the distributed setting of Algorithm~\ref{alg:ADMM-ForOurProblem} using the ADMM method, it also provides the potential benefits associated with federated learning and distributed systems~\cite{li2020federated, lowy2023private}. 

\vspace{0.1cm}

In the relaxed version of problem~\eqref{eq:optModel1}, different solutions $\bS_i^{\ast}$ with a fixed $\bS_i^{\ast} \mathbf{1}=\bu^{\ast}$ yield the same objective value if $\bS_i^{\ast}$ satisfies all the constraints. Thus, potentially infinitely many solutions to our convex problem exist, and many are not binary. To promote a binary solution for the final decision, we introduce the following regularizer into the objective function of the relaxed version of problem~\eqref{eq:optModel1}:
\begin{equation}
\begin{split}
    \Re(\bS) = - \frac{\tilde{\lambda}}{2} \sum_{r=1}^{\cR} \sum_{t=1}^{|\bT|} \sum_{i=1}^{n}(\bS_i)_{r, t}((\bS_i)_{r, t} - 1)
\end{split}
\normalsize
\end{equation}
\normalsize
where $\tilde{\lambda} \in \bbR_{+}$ is the regularization parameter and $(\bS_i)_{(r, t)} \in [0, 1]$. 
This regularizer has the effect of driving the elements of matrix $\bS$ towards the binary domain $\{0, 1\}$. The regularizer penalizes any deviations from this domain in the objective function. While convexity is sacrificed due to regularization, ADMM can still be convergent in nonconvex problems~\cite{hong2016convergence}.

\vspace{0.1cm}

Algorithm~\ref{alg:ADMM-ForOurProblem} solves the relaxed version of problem~\eqref{eq:optModel1}. Since the solution to the relaxed version of problem~\eqref{eq:optModel1} may not be binary (due to relaxation), we need to project it back to the feasible region. For computational purposes, we suggest using $\ell_1$ projection by solving the following mixed integer (linear) problem (MILP)
\begin{equation}   
\label{eq:linearModelS}
    \begin{split}
        \min_{\{\bS_i, c_i\}_{i=1}^{n}} \quad & \sum_{i=1}^{n}\| \bS_{i}\mathbf{1} - \bu_i^{\ast} \|_{1} \\
        \stt \quad & \bD \bS_i \mathbf{1} = \bq_i,\quad \forall i=1, 2, \dots, n \\
        & \bS_i^{\top} \mathbf{1} = \mathbf{1},\quad  \forall i=1, 2, \dots, n \\
        & \bS_i \in \lbrace0, 1\rbrace^{(|\cR|\cdot|\bT|)\times(|\cN_{i}|)},\quad  \forall i=1, 2, \dots, n \\
        & \bS_i^{\top} \bdelta \leq \bb_i \odot \bB_i \boldsymbol{\eta},\quad  \forall i=1, 2, \dots, n\\
        & c_i \geq \alpha_i(\bdelta^{\top} \bS_i \mathbf{1} - \gamma_i),\quad \forall i=1, 2, \dots, n \\
        & c_1 + c_2 + \dots + c_n \leq \Omega \\
        & c_i \geq 0,\quad \forall i=1, 2, \dots, n \\
    \end{split}
\normalsize
\end{equation}
where $\bu_i^{\ast}, \forall i=1,2,\dots,n$ is the optimal solution obtained by Algorithm~\ref{alg:ADMM-ForOurProblem}. 
Clearly, this problem can be reformulated as a MILP problem and solved using off-the-shelf solvers like Gurobi. Solving problem~\eqref{eq:linearModelS} can be easier than  problem~\eqref{eq:optModel1}. Problems~\eqref{eq:optModel1} and~\eqref{eq:linearModelS} have the same variable size and similar constraints, but the objective functions are different. 
While the objective function in problem~\eqref{eq:linearModelS} can be restructured as a linear programming problem, we have a polynomial objective function in problem~\eqref{eq:optModel1} that introduces more complexity. 

\vspace{0.2cm}

\section{Experiments} \label{sec:NumericalExperiments}
We evaluate our incentive scheme's effectiveness using Los~Angeles area data. The presence of multiple routes between most origin-destination (OD) pairs makes the Los~Angeles area particularly suitable for our assessment. We use the data collected by the Archived Data Management System (ADMS), a comprehensive transportation dataset compilation by University of Southern California researchers~\cite{anastasiou2019admsv2}. This system aggregates data from Los~Angeles, Orange, San Bernardino, Riverside, and Ventura Counties, offering a robust data source for analysis.

For our evaluations, we need to estimate the OD matrix. The $(i, j)$-th entry of the OD matrix represents the count of drivers traveling between origin $i$ and destination $j$. We need to estimate the OD matrix using the available network flow information due to the unavailability of drivers' routing data. 
The OD matrix estimation problem is challenging due to its under-determined nature~\cite{van1980most, cascetta1984estimation, bell1991real}. OD matrices are categorized as static or dynamic~\cite{bera2011estimation}. However, many dynamic OD estimation (DODE) methods are computationally impractical for our high-resolution data. Additionally, some studies rely on existing OD matrix data~\cite{krishnakumari2019data, carrese2017dynamic, nigro2018exploiting, kim2014using}, which we lack access to. Given these constraints, we adopt the OD estimation algorithm proposed by \cite{ma2018estimating}.
All the codes are publicly available at: \href{https://github.com/ghafeleb/Incentive_Systems_for_New_Mobility_Services}{https://github.com/ghafeleb/Incentive\_Systems\_for\_New\-\_Mobility\_Services}.

\subsection{Simulation Model}
\label{sec:simModel}
\begin{figure}[] 
  \centering
  \begin{tabular}{@{}c@{}}
    \includegraphics[width=0.4\linewidth]{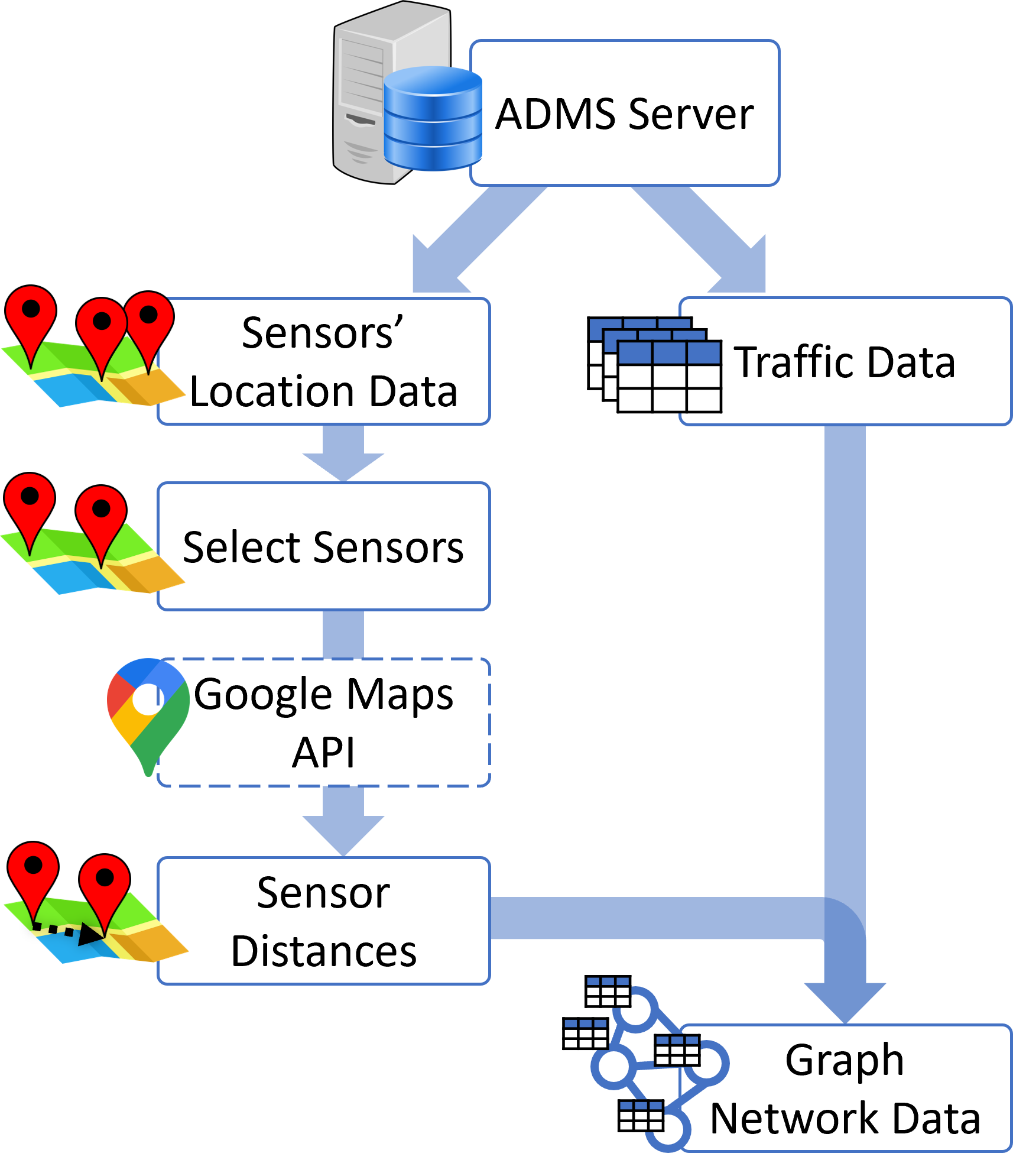} \\[\abovecaptionskip]
  \end{tabular}
  \vspace{-0.45cm}
  \caption{Data preparation workflow.}
  \label{fig:data_preparation_workflow}
\end{figure}

\begin{figure}[] 
  \centering
  \begin{tabular}{@{}c@{}}
    \includegraphics[width=0.75\linewidth]{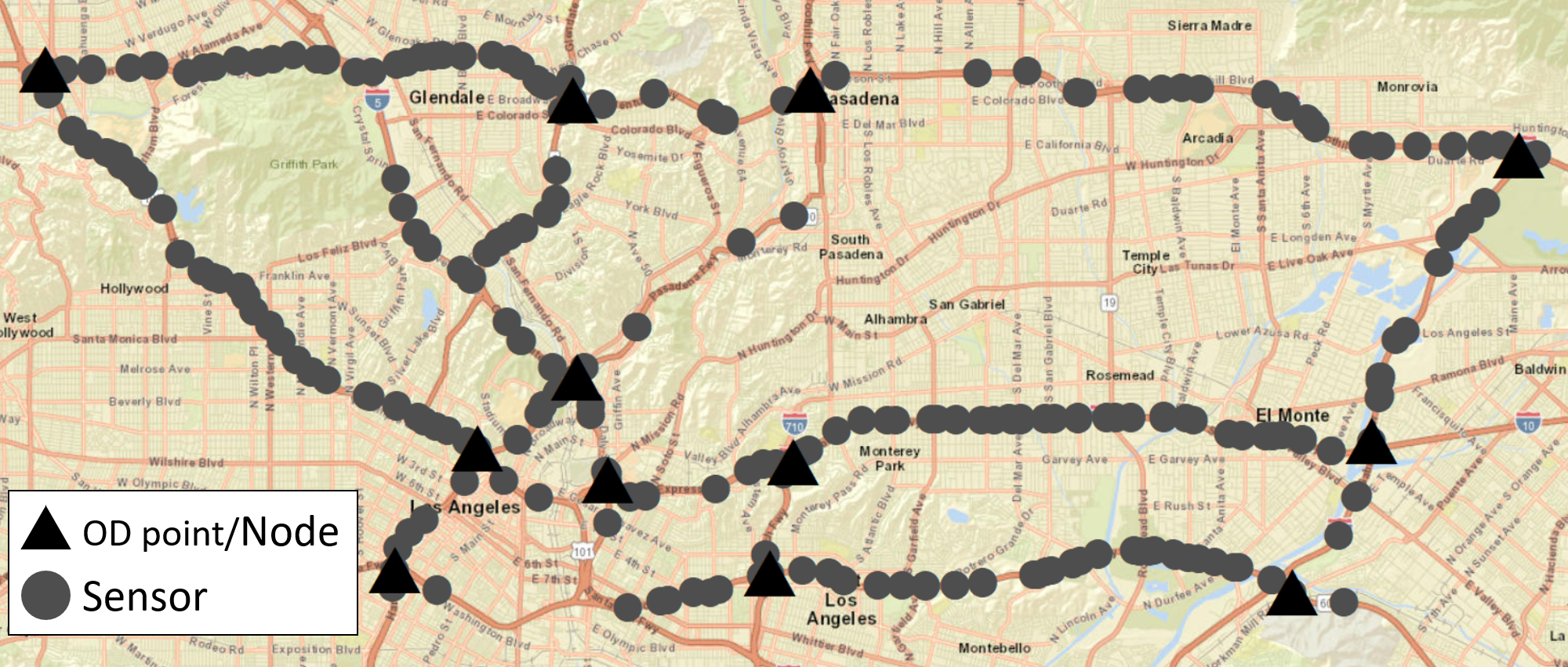} \\[\abovecaptionskip]
  \end{tabular}
  \vspace{-0.5cm}
  \caption{Studied region and the highway sensors inside the region.}
  \label{fig:region_y3}
\end{figure}

\begin{figure}[] 
  \centering
  \begin{tabular}{@{}c@{}}
    \includegraphics[width=0.45\linewidth]{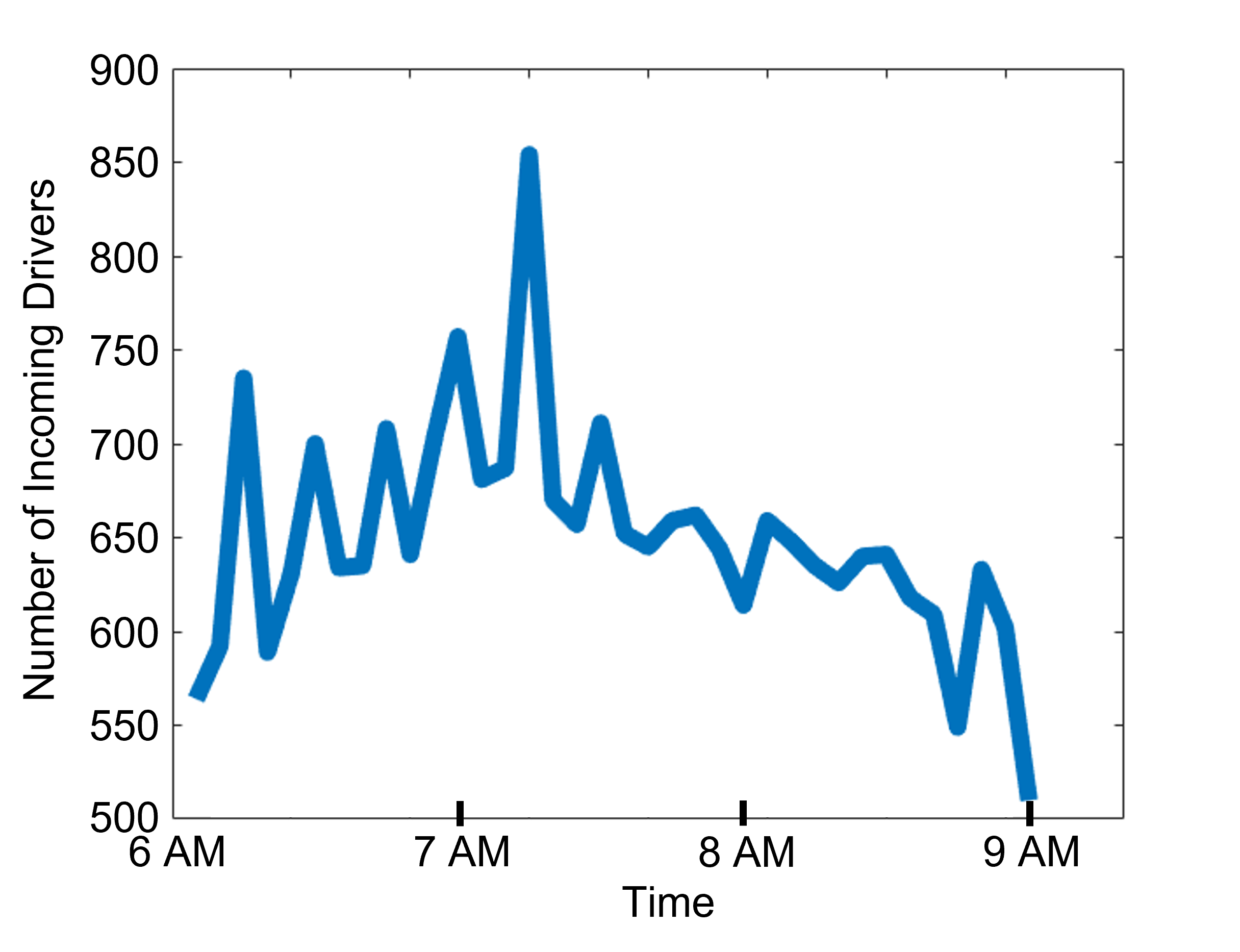} \\[\abovecaptionskip]
  \end{tabular}
  \vspace{-0.6cm}
  \caption{Total estimated number of drivers entering the system over time (in 5-minute intervals).}
  \label{fig:number_driver}
\end{figure}

First, we extract sensor details, including their locations. We extract the speed and volume data of selected sensors. Nodes for the network graph are chosen from on-ramps and highway intersections. Connecting link data is derived from in-between sensors. Node distances are determined via Google Maps API. Data preparation workflow is shown in Figure~\ref{fig:data_preparation_workflow}.
The network under consideration includes highways surrounding Downtown Los Angeles, as depicted in Figure~\ref{fig:region_y3}, and consists of 12 nodes, 32 links, and a total road length of 288.1 miles.
We have 144 OD pairs, and we employ the algorithm from~\cite{ma2018estimating} on the network's speed and volume data to estimate OD pairs. 
Figure~\ref{fig:number_driver} shows the total estimated incoming drivers per time interval. 
We explore 3 routing options for each OD pair. Initially, the shortest path is determined. Subsequently, links in the first path are removed to uncover the second shortest path if available. This process is repeated for the third route. Based on this process, we find 270 paths between OD pairs.

\vspace{0.2cm}

We focus on incentivizing the organizations to change their behavior for the 7~AM to 8~AM interval (which is the rush hour based on the estimated number of incoming drivers in Figure~\ref{fig:number_driver}). Although we have selected 7~AM to 8~AM as the incentivization time period, we also include 8~AM to 8:30~AM in our experiments because some of the drivers entering between 7~AM and 8~AM may not finish their route before 8~AM. To track the effect of these drivers on the total travel time of the system, we include traffic flow from 8~AM to 8:30~AM in our analysis as well. 
The OD estimation algorithm's projected total count of drivers entering the system from 6 AM to 9 AM is illustrated in Figure~\ref{fig:number_driver}.
From 7 AM to 8:30 AM, a total of 11985 drivers enter the system.

\vspace{0.2cm}

We consider the traffic volume of the network at UE in our baseline.
To compute the volume of the network at UE, we use the UE~algorithm in~\cite{ghafelebashi2023congestion}. The algorithm receives the volume (historical data) and OD estimation as inputs and returns the matrix $\bR_{\text{UE}}$ and route travel time vector $\bdelta_{\text{UE}}$ at UE.
To compute the cost of organizations' incentivization, we need to know the route travel times when drivers have made decisions based on the UE route travel time $\bdelta_{\text{UE}}$. Hence, we compute the new volume vector $\bv_{\text{new}} = \bR_{\text{UE}} \bS_{\text{UE}} \mathbf{1}$ where $\bS_{\text{UE}}$ is the assignment of drivers to the fastest route based on the UE route travel time vector $\bdelta_{\text{UE}}$. Using the BPR function, volume vector $\bv_{\text{new}}$, and $\bdelta_{\text{UE}}$, we compute $\bdelta$ that denotes the travel time of the routes if drivers make decision based on $\bdelta_{\text{UE}}$ and $\boldsymbol{\eta}$ denotes the minimum travel time between the different OD pairs.

\subsection{Results}
In this subsection, using our model and algorithm, we study the impact of organization incentivization for different budget values, the number of organizations, VOTs, and the percentage of drivers who are employed by the organizations in the incentivization program. The remaining drivers are assumed to be background drivers who follow the $\bdelta_{\text{UE}}$. 
We consider two scenarios for the percentage of drivers:
\begin{itemize}
    \item Scenario~I: Among the drivers entering the system between 7~AM and 8~AM, 10\% of them (i.e., 812 drivers) belong to organizations that we can incentivize. 
    \item Scenario~II: Among the drivers entering the system between 7~AM and 8~AM, 20\% of them (i.e., 1624 drivers) belong to organizations that we can incentivize.  
\end{itemize} 
Drivers in each organization are selected uniformly at random, and all selected drivers of Scenario~I are included in Scenario~II to have a fair comparison between the two scenarios.
Our base VOT is derived from the estimation of \cite{castillo2020benefits}, which is \$2.63 per minute or \$157.8 per hour. The default number of organizations in our experiments is 10. 

\vspace{0.2cm}

\begin{figure}[] 
  \centering
  \begin{tabular}{@{}c@{}}
    \includegraphics[width=0.45\linewidth]{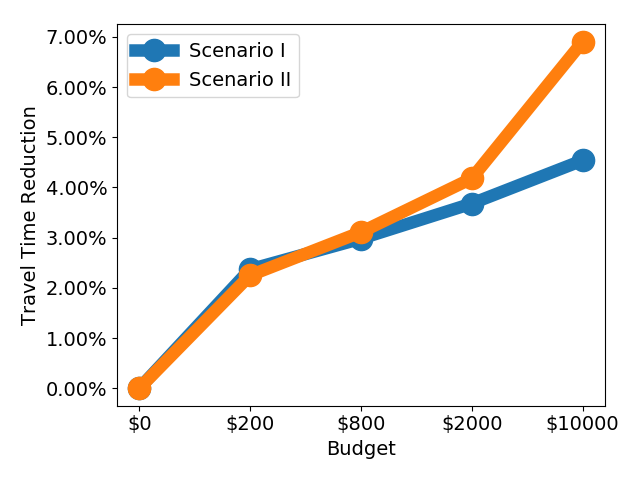} \\[\abovecaptionskip]
  \end{tabular}
  \vspace{-0.6cm}
  \caption{Percentage of travel time decrease with different budgets at VOT=\$157.8 per hour.}
  \label{fig:TTDecrease2.63}
\end{figure}

\begin{figure}[] 
  \centering
  \begin{tabular}{@{}c@{}}
    \includegraphics[width=0.45\linewidth]{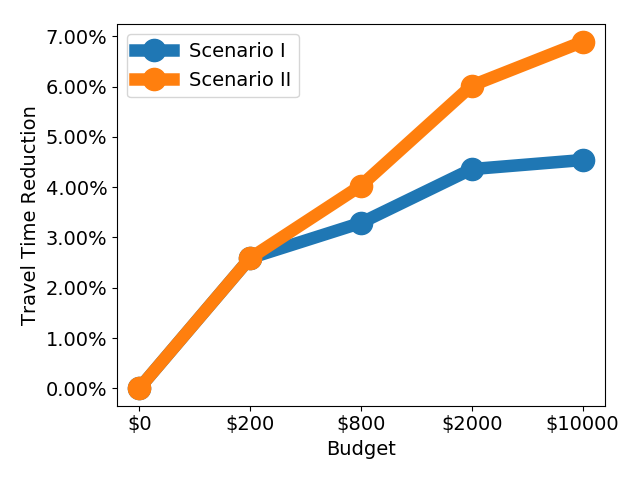} \\[\abovecaptionskip]
  \end{tabular}
  \vspace{-.6cm}
  \caption{Percentage of travel time decrease with different budgets at VOT=\$78.9 per hour.}
  \label{fig:TTDecrease1.315}
\end{figure}

The percentage of travel time decrease with incentivization as compared to a system with no incentivization scheme with VOT of $\$157.8$ is presented in Figure~\ref{fig:TTDecrease2.63} for Scenario~I and Scenario~II. In our plots, the budget of \$0 shows the case of a no-incentivization platform. 
The no-incentivization system solution essentially assumes all drivers are background drivers. We observe that by increasing the available budget, the amount of decrease in travel time increases (as expected). This decrease is more in larger budgets in Scenario~II because the model has access to more drivers to select and has more flexibility to recommend alternative routes. 
For the purpose of sensitivity analysis, we also provide travel time decrease for both Scenario~I and Scenario~II with a different VOT of $\frac{\$157.8}{2} = \$78.9$ per hour in Figure~\ref{fig:TTDecrease1.315}. The comparison of results for different VOTs in Figure~\ref{fig:TTDecrease2.63} and Figure~\ref{fig:TTDecrease1.315} shows that for a very large budget, the decrease in travel time is almost similar. This is because none of the models utilize the entire budget at a $\$10,000$ budget. However, with a smaller VOT and a budget of \$2000, there is a large gap between  Scenario~I and Scenario~II  because Scenario~II has access to more drivers to deviate. 

\begin{figure}[] 
  \centering
  \begin{tabular}{@{}c@{}}
    \includegraphics[width=0.45\linewidth]{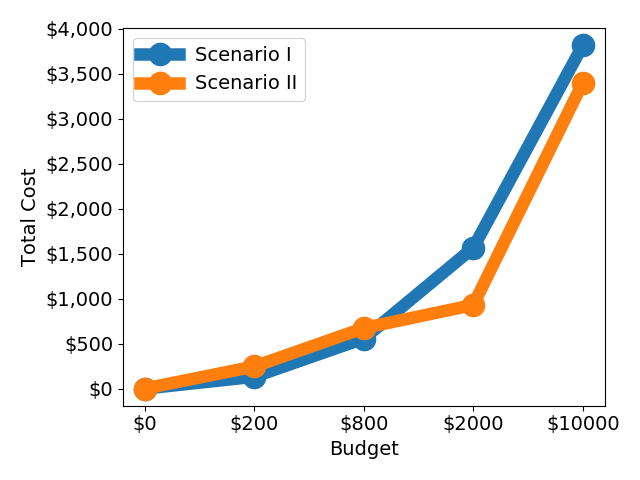} \\[\abovecaptionskip]
  \end{tabular}
  \vspace{-0.6cm}
  \caption{Total cost of incentivization of 10 organizations with different budgets in Scenario~I and Scenario~II and VOT=\$157.8 per hour.}
  \label{fig:cost157.8}
\end{figure}

\begin{figure}[] 
  \centering
  \begin{tabular}{@{}c@{}}
    \includegraphics[width=0.45\linewidth]{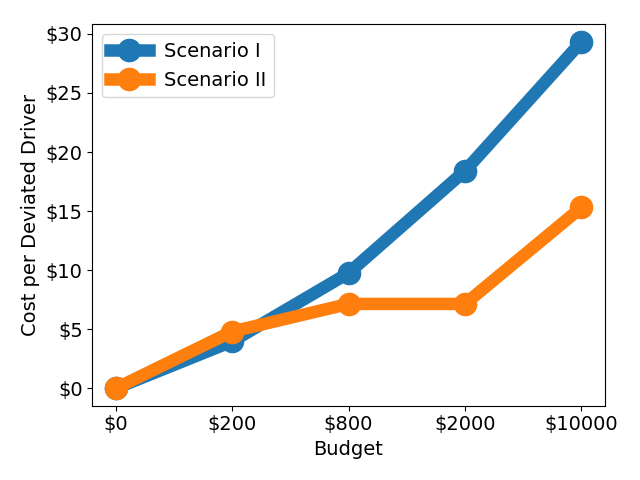} \\[\abovecaptionskip]
  \end{tabular}
  \vspace{-0.6cm}
  \caption{Cost of incentivization per deviated drivers of 10 organizations with different budgets in Scenario~I and II and VOT=\$157.8 per hour.}
  \label{fig:costPerDriver157}
\end{figure}

For the next analyses of our numerical results, we only report the results for our base VOT (\$2.63 per minute or \$157.8 per hour) because the results follow similar patterns with VOT of \$78.9. 
In Figure~\ref{fig:cost157.8}, we present the total incentivization cost for different budgets in both Scenario~I and Scenario~II when there are 10 organizations in the system. This cost increases when the available budget is more. This pattern shows that the platform can utilize the resources when it has access to more money. We observe that more involvement of drivers leads to a slightly smaller cost because in Scenario~II model has more flexibility in selecting drivers.
Figure~\ref{fig:costPerDriver157} shows the cost per deviated driver for the two scenarios. Although the gap between the total cost of Scenario~I and Scenario~II is small, the cost per driver is significantly smaller in Scenario~II because the model has more flexibility in choosing the drivers efficiently. 
Moreover, the cost per driver increases with the budget. This shows that our model utilizes our budget efficiently by providing more affordable incentives first when the budget is low.
As Table~\ref{table:nDeviatedDriver} shows, the number of incentivized drivers in Scenario~II is larger because there are more drivers for selection. 
\begin{table}[]
\centering
\begin{tabular}{|c|cccc|}
\hline
\multirow{2}{*}{Scenario} & \multicolumn{4}{c|}{Budget}                                                                     \\ \cline{2-5} 
                          & \multicolumn{1}{c|}{\$200} & \multicolumn{1}{c|}{\$800} & \multicolumn{1}{l|}{\$2000} & \$10000 \\ \hline
I                         & \multicolumn{1}{c|}{33}    & \multicolumn{1}{c|}{57}    & \multicolumn{1}{c|}{85}     & 130     \\ \hline
II                        & \multicolumn{1}{c|}{51}    & \multicolumn{1}{c|}{94}    & \multicolumn{1}{c|}{130}    & 222     \\ \hline
\end{tabular}
\vspace{0.4cm}
\caption{Distribution of the number of drivers that were assigned to an alternative route.}
\label{table:nDeviatedDriver}
\end{table}

\begin{figure}[H] 
  \centering
  \begin{tabular}{@{}c@{}}
    \includegraphics[width=0.45\linewidth]{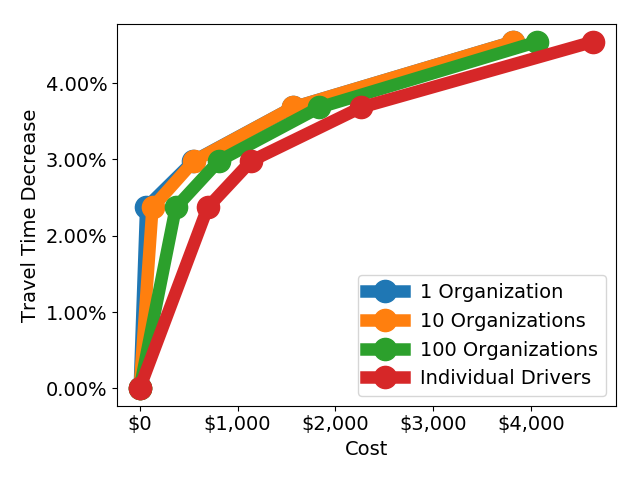} \\[\abovecaptionskip]
  \end{tabular}
  \vspace{-0.6cm}
  \caption{Travel time decrease vs. incentivization cost for different number of organizations in Scenario~I and VOT=\$157.8 per hour.}
  \label{fig:costTT_ScenarioI_157}
\end{figure}

\begin{figure}[H] 
  \centering
  \begin{tabular}{@{}c@{}}
    \includegraphics[width=0.45\linewidth]{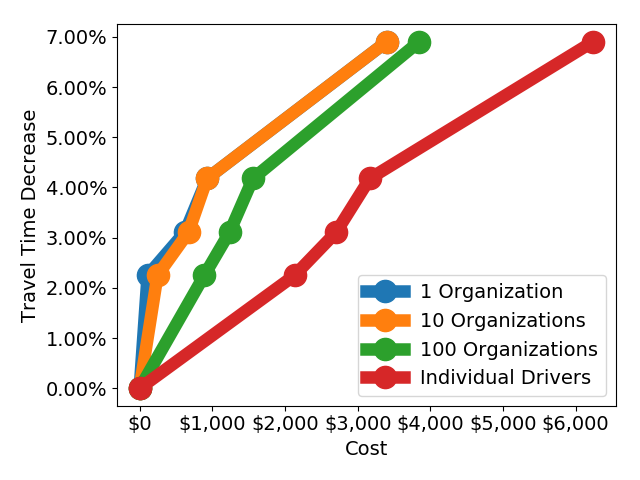} \\[\abovecaptionskip]
  \end{tabular}
  \vspace{-0.6cm}
  \caption{Travel time decrease vs. incentivization cost for different number of organizations in Scenario~II and VOT=\$157.8 per hour.}
  \label{fig:costTT_ScenarioII_157}
\end{figure}

The number of organizations in the system can alter the total travel time and cost. Figure~\ref{fig:costTT_ScenarioI_157} and Figure~\ref{fig:costTT_ScenarioII_157} illustrate the percentage decrease of travel time and total cost when there are different number of organizations in the system.  As an extreme case, we also include the case that each organization contains one driver (i.e., we incentivize individuals rather than organizations). In Figure~\ref{fig:costTT_ScenarioI_157} and Figure~\ref{fig:costTT_ScenarioII_157}, we observe a larger cost for reducing the same amount of travel time decrease when there are more organizations in the system. The intuitive reason behind this observation is as follows. For each organization, after incentivization, some drivers lose time, and some gain travel time. At the organizational level, the time changes of drivers can cancel each other out, and hence we may not need to compensate the organization significantly. When the number of drivers per organization decreases, the canceling effect becomes weaker, and the incentivization costs more. This is in line with our discussion in Section~\ref{sec:whyOrganizations}. This also explains why incentivizing organizations is much more cost-efficient than incentivizing individual drivers.

\section{Conclusion} 
\label{sec:conclusion}
In this paper, we study the problem of incentivizing organizations to reduce traffic congestion. To this end, we developed a mathematical model and provided an algorithm for offering organization-level incentives. In our framework, a central planner collects the origin-destination and routing information of the organizations. Then, the central planner utilizes this information to offer incentive packages to organizations to incentivize a system-level optimal routing strategy. 
In particular, we focused on minimizing the total travel time of the network. However, other utilities can be used in our framework. Finally, we employed data from the Archived Data Management System (ADMS) to evaluate the performance of our model and algorithm in a representative traffic scenario in the Los~Angeles area. A 6.90\% reduction in the total travel time of the network was reached by our framework in the experiments. More importantly, we observed that incentivizing companies/organizations is more cost-efficient than incentivizing individual drivers. 
As future work, it is important to study the effect of incentivization to change the start time of the trip.  This is particularly relevant in future mobility services because many of them, such as delivery services, are flexible in terms of trip time to a certain degree.
In addition, we can consider the stochastic nature of making decisions in routing by individual drivers. Moreover, we can extend the incentivization framework to the case that not all organizations accept their received offer. 

\bibliographystyle{unsrtnat}
\bibliography{arxiv}

\appendix

\section{List of Notations} \label{sec:notation}
\label{apdx:notations}
Traffic network spatiotemporal parameters:
\begin{itemize}
    \item $\cG$: Directed graph of the traffic network
    \item $\cV$: Set of nodes of graph $\cG$ which correspond to major intersections and ramps
    \item $\cE$: Set of edges of graph $\cG$ which correspond to the set of road segments
    \item $|\mathcal{E}|$: Total number of road segments/edges in the network $\cG$ (i.e. the cardinality of the set $\mathcal{E}$)
    \item $\ell$: An edge of graph $\cG$ which corresponds to a link/road segment in the traffic network
    \item $\mathcal{R}_j$: Set of possible route options for driver~$j$
    \item $\mathcal{R}$: Total set of possible route options for all OD pairs
    \item $|\mathcal{R}|$: Total number of possible route options (i.e. the cardinality of the set $\mathcal{R}$)
    \item $\mathbf{r}$: Route vector 
    \item $\mathbf{T}$: Set of time of periods
    \item $|\mathbf{T}|$: Number of time units (i.e. the cardinality of $\mathbf{T}$)
    \item $\theta_{\ell, t}$: Travel time of link $\ell$ at time $t$
    \item $F(.)$: Total travel time function
    \item $T_{\br}$: The travel time for route $\br$ 
\end{itemize}

\vspace{0.2cm}

BPR function and its parameters:
\begin{itemize}
    \item $f_{\text{BPR}}(.)$: BPR function
    \item $v$: The traffic volume of the link
    \item $w$: The practical capacity of the link
    \item $\theta_0$: The free flow travel time of the link
\end{itemize}

\vspace{0.2cm}

Optimization model parameters:
\begin{itemize}
    \item $\mathcal{N}_i$: Set of drivers of organization~$i$
    \item $|\mathcal{N}_i|$: Total number of drivers of organization~$i$ (i.e. the cardinality of set $\mathcal{N}_i$)
    \item $\mathcal{N}$: Set of all drivers
    \item $|\mathcal{N}|$: Total number of drivers (i.e., the cardinality of set $\mathcal{N}$)
    \item $\mathbf{v}_t$: Volume vector of road segments at time $t$ 
    \item $\hat{\bv}$: The vector of the estimated volume of links at different times in the horizon
    \item $\hat{v}_{\ell, t}$: The $(|\cE|\times t + \ell)^{th}$ element of vector $\hat{\bv}$ representing the volume of the $\ell^{th}$ link at time $t$
    \item $\bR$: The matrix of the probability of a driver being at each link given their route
    \item $\br_{\ell, t}$: The row of matrix $\bR$ that corresponds to link $\ell$  at time $t$
    \item $\bD$: The matrix of route assignments of the OD pairs
    \item $\bq_i$: The vector of number of drivers of organization~$i$ for each OD pair
    \item $\bdelta$: The vector of travel time of routes at different times
    \item $\boldsymbol{\eta}$: The vector of shortest travel time between different OD pairs at different times
    \item $\bb_i$: This vector contains the factors by which the travel time of assigned routes can be larger than the shortest travel time of the drivers of organization~$i$
    \item $\bB_i$: The matrix of shortest travel time assignment of drivers of organization~$i$
    \item $\alpha_i$: VOT for organization~$i$
    \item $\balpha$: The vector of VOT values for the different organizations
    \item $\gamma_i$: Total travel time of organization~$i$ in the absence of incentivization platform
    \item $\Omega$: Budget for incentivization
    \item $K$: The number of OD pairs
\end{itemize}

\vspace{0.2cm}

Decision variables:
\begin{itemize}
    \item $s_i^{\mathbf{r},j}$: Decision parameter indicates whether route~$\br$ is assigned to driver~$j$ from organization~$i$
    \item $\bs_i^j$: The binary route assignment vector of driver~$j$ from organization~$i$ 
    \item $\bS_i$: Decision matrix of drivers of organization~$i$
    \item $\bS$: Decision matrix of all drivers
    \item $c_i$: The cost of incentive offered to organization~$i$
\end{itemize}

\section{Model and Notations: An Illustrative Example}
\label{appdx:ModelNotationExample}
\begin{figure}[H]
  \centering
  \begin{tabular}{@{}c@{}}
    \includegraphics[width=.2\linewidth]{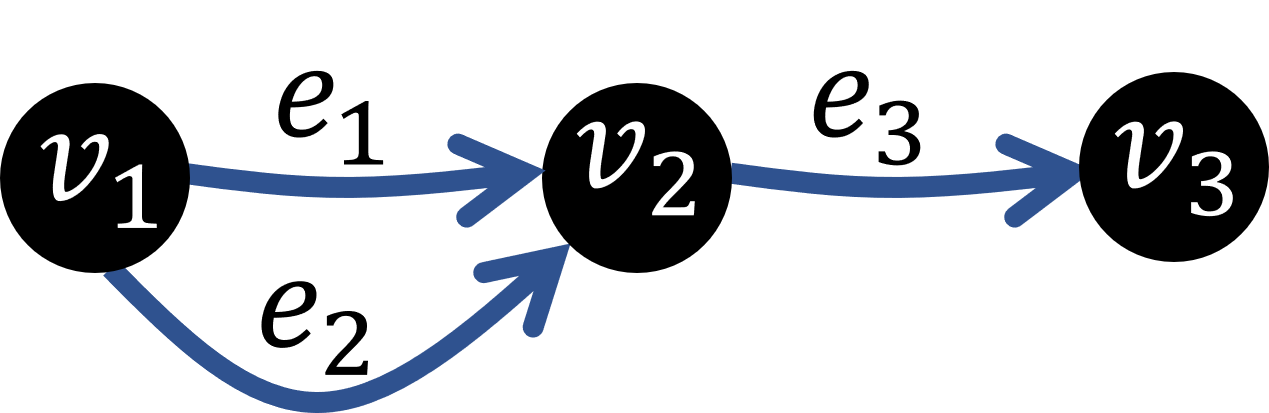} \\[\abovecaptionskip]
  \end{tabular}
  \caption{Network example $\cG_1$.}\label{fig:network}
\end{figure}

\normalsize
\noindent{\normalsize This section presents an example of a network to elaborate on our model and its associated notations.  Consider  network $\cG_1$ in Figure~\ref{fig:network}}
\vspace{-.2cm}
\noindent{\normalsize where $\cV = \{\nu_1, \nu_2, \nu_3\}$ is the set of nodes and $\cE = \{e_1, e_2, e_3\}$ is the set of edges (roads). The link details and their attributes are shown in Table~\ref{fig:tableLinks}. The origin-destination (OD) pair is ($\nu_1$, $\nu_3$). There are two routes from the origin to the destination, detailed in Table~\ref{fig:tableRoutes}. $\bT = \{1, 2, 3\}$ is the set of time of periods.  Each period is 0.2 hours. To estimate driver locations at each time,  we need matrix $\bR \in [0, 1]^{9\times6}$ as follows}

\begin{figure}[H]
\begin{align*}
\bR = 
\begin{blockarray}{ccccccc}
& \makecell{t_1\\=1\\ \br_1} & \makecell{t_1\\=1\\ \br_2} & \makecell{t_1\\=2\\ \br_1} & \makecell{t_1\\=2\\ \br_2} & \makecell{t_1\\=3\\ \br_1} & \makecell{t_1\\=3\\ \br_2} \\
\begin{block}{c(cccccc)}
  \makecell{t_2=1, e_1} & 1 & 1 & 0 & 0 & 0 & 0 \\
  \makecell{t_2=1, e_2} & 0 & 0 & 0 & 0 & 0 & 0 \\
  \makecell{t_2=1, e_3} & 0.5 & 0 & 0 & 0 & 0 & 0 \\
  \makecell{t_2=2, e_1} & 0 & 0 & 1 & 1 & 0 & 0 \\
  \makecell{t_2=2, e_2} & 0 & 1 & 0 & 0 & 0 & 0 \\
  \makecell{t_2=2, e_3} & 0.5 & 0 & 0.5 & 0 & 0 & 0 \\
  \makecell{t_2=3, e_1} & 0 & 0 & 0 & 0 & 1 & 1 \\
  \makecell{t_2=3, e_2} & 0 & 0 & 0 & 1 & 0 & 0 \\
  \makecell{t_2=3, e_3} & 0 & 0 & 0.5 & 0 & 0.5 & 0 \\
\end{block}
\end{blockarray}
\end{align*}
\end{figure}
\normalsize
\noindent
where $t_1$ represents the time when the driver enters the system and $t_2$ is the time when the driver reaches the road.
\begin{table}[H]
\centering
\begin{tabular}{c|c|c|c|}
\cline{2-4}
                       &  \makecell{Length\\ (Mile)} & \makecell{Speed\\ (mph)} & \makecell{Travel Time\\ (Hour)} \\ \hline
\multicolumn{1}{|c|}{$e_1$} &  5 & 50 & 0.1\\ \hline
\multicolumn{1}{|c|}{$e_2$} &  10 & 50 & 0.2\\ \hline
\multicolumn{1}{|c|}{$e_3$} &  5 & 50 & 0.1\\ \hline
\end{tabular}
\vspace{+0.35cm}
\caption{Set of edges.} \label{fig:tableLinks}
\end{table}
\normalsize

\begin{table}[H]
\centering
\begin{tabular}{c|c|c|}
\cline{2-3}
                       &  $\br$ & Graph\\ \hline
\multicolumn{1}{|c|}{\makecell{Route 1\\ $e_1 \rightarrow e_3$}} &  $\br_1=\begin{bmatrix} 1 \\ 0 \\ 1\end{bmatrix}$ & \includegraphics[width=.2\linewidth]{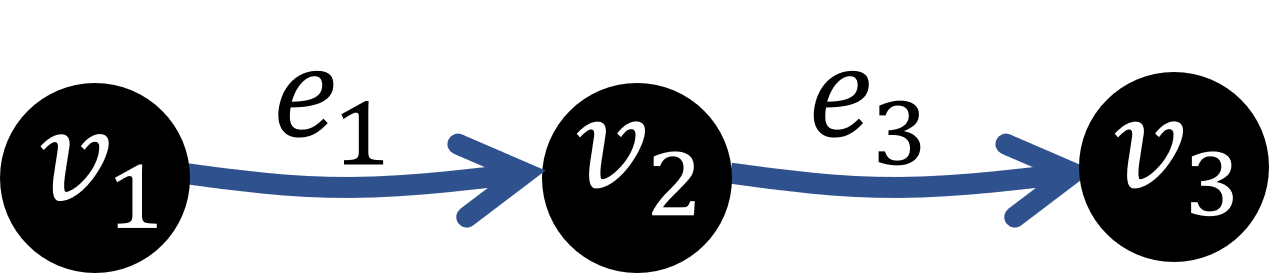} \\ \hline
\multicolumn{1}{|c|}{\makecell{Route 2\\ $e_2 \rightarrow e_3$}} & $\br_2 =\begin{bmatrix} 0 \\ 1 \\ 1\end{bmatrix}$ & \includegraphics[width=.2\linewidth]{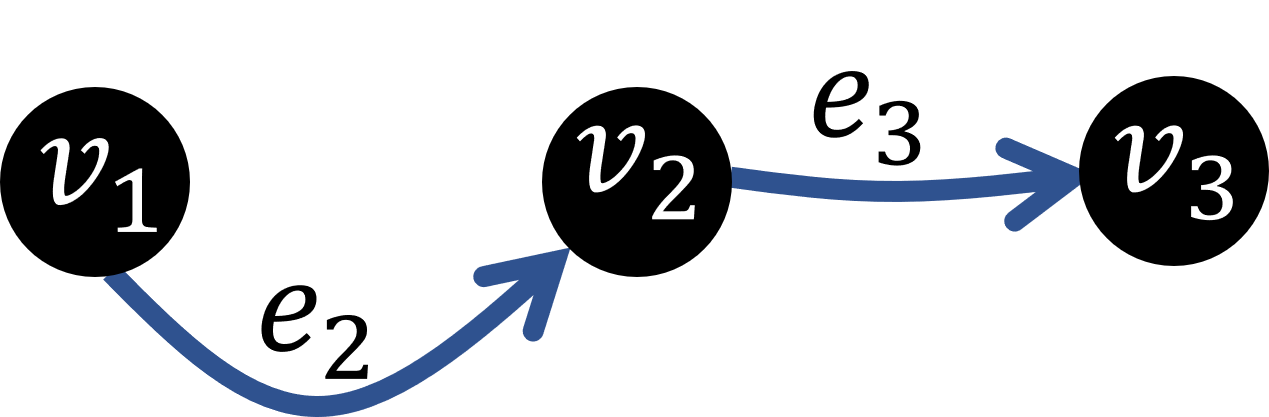} \\ \hline
\end{tabular}
\vspace{+0.35cm}
\caption{Set of routes.} \label{fig:tableRoutes}
\end{table}
\normalsize
\vspace{-0.8cm}
Assume there are two organizations in the system. VOT for organization~1 is \$2.50/minute and VOT for organization~2 is \$1.50/minute so $\alpha_1=2.50$ and $\alpha_2=1.50$. There are three drivers in the system and $\cN = \{d_1, d_2, d_3\}$ such that $\cN_1 = \{d_1, d_2\}$ and $\cN_2 = \{d_3\}$. Drivers $d_1$ and $d_2$ start their travel at time 1 and the travel of driver $d_3$ starts  at time 2.
Given that drivers $d_1$ and $d_2$ have to finish their travel as soon as possible and the travel time of driver $d_3$ can be delayed up to two times the shortest travel time, matrices $\bB_1$ and $\bB_2$ and vectors $\bb_1$ and $\bb_2$ will be defined as follows \vspace{-0.5cm}
\noindent
\begin{figure}[H]
  \centering
  \begin{minipage}[b]{0.4\textwidth}
    \begin{align*}
    \bB_1 = 
    \begin{blockarray}{cccc}
    & \multicolumn{3}{c}{\;\;\text{OD Assignment}} \\
    \begin{block}{c(ccc)}
      \makecell{d_1} & 1 & 0 & 0  \\
      \makecell{d_2} & 1 & 0 & 0  \\
    \end{block}
    \end{blockarray}
    \end{align*}
  \end{minipage}%
  \begin{minipage}[b]{0.4\textwidth}
    \begin{align*}
    \bb_1 = 
    \begin{blockarray}{cc}
    & \makecell{\text{Travel Time Multiplier}}  \\
    \begin{block}{c(c)}
      \makecell{d_1} & 1  \\
      \makecell{d_2} & 1  \\
    \end{block}
    \end{blockarray}
    \end{align*}
  \end{minipage}
\end{figure}
\vspace{-1.0cm}
\begin{figure}[H]
  \centering
  \begin{minipage}[b]{0.4\textwidth}
    \begin{align*}
    \bB_2 = 
    \begin{blockarray}{cccc}
    & \multicolumn{3}{c}{\;\;\text{OD Assignment}}  \\
    \begin{block}{c(ccc)}
      \makecell{d_3} & 0 & 1 & 0  \\
    \end{block}
    \end{blockarray}
    \end{align*}
  \end{minipage}%
  \begin{minipage}[b]{0.4\textwidth}
    \begin{align*}
    \bb_2 = 
    \begin{blockarray}{cc}
    & \makecell{\text{Travel Time Multiplier}}  \\
    \begin{block}{c(c)}
      \makecell{d_3} & 2  \\
    \end{block}
    \end{blockarray}
    \end{align*}
  \end{minipage}
\end{figure}

\vspace{-0.8cm}

\section{Distributed Incentivization Algorithm}\label{appdx:alg}
\begin{algorithm}[] 
    \caption{Distributed Organization-Level Incentivization via ADMM}
    \label{alg:ADMM-ForOurProblem}
	\begin{algorithmic}[1]
	 \State \textbf{Input}: Initial values: $\bomega^0$, $\bS_i^0$, $\bH_i^0$, $\bW_i^0$, $\bZ_i^0$, $\bu^0$, $\beta_i^{0}$, $\tilde{\beta}^{0}$, $\tilde{\bc}^{0}$, $\bla^0_{1, i} \in \mathbb{R}^{|\cR|\cdot|\bT|}$,  $\bla^0_2 \in \mathbb{R}^{|\cE| \cdot |\bT|}$,  $\bla^0_{3, i} \in \mathbb{R}^{K \cdot |\bT|}$,  $\bla^0_{4, i} \in \mathbb{R}^{|\cR| \cdot |\bT|}$,  $\bLa^0_{5, i} \in \mathbb{R}^{(|\cR|\cdot|\bT|) \times |\cN_i|}$,  $\bla^0_{6, i} \in \mathbb{R}^{|\cN_i|}$,  $\bla^0_7 \in \mathbb{R}^{n}$,   $\bLa^0_{8, i} \in \mathbb{R}^{(|\cR|\cdot|\bT|) \times |\cN_i|}$,   $\lambda^0_9 \in \mathbb{R}$,   $\bLa^0_{10, i} \in \mathbb{R}^{(|\cR|\cdot|\bT|) \times |\cN_i|}$, Dual update step: $\rho$, Number of iterations: $\tilde{T}$.
        \For {$t = 0, 1, \ldots, \tilde{T}$}
	    	    
	    \For {$\ell = 0, 1, \ldots, |\cE|$}
	    \For {$\hat{t} = 1, \ldots, |\bT|$}
	    \State $\omega^{t+1}_{\ell, \hat{t}} = \argmin\limits_{\omega_{\ell, \hat{t}}}   \omega_{\ell, \hat{t}} \theta_{\ell, t}(\omega_{\ell, \hat{t}})  + \bla^{t}_{2, (\ell, \hat{t})}(\omega_{\ell, \hat{t}} - \br_{\ell, \hat{t}} \left(\sum_{i=1}^{n}\bu_i^{t}\right)) + \frac{\rho}{2} (\omega_{\ell, \hat{t}} - \bR_{\ell, \hat{t}} \left(\sum_{i=1}^{n}\bu_i^{t}\right))^{2}$ \label{step:omega}
	    \EndFor
	    \EndFor
	    
	    \For {$i = 1, \ldots, n$}
	    \State $\bS^{t+1}_i = (-\bla_{1,i}^{t}\mathbf{1}^{\top} - \bLa_{5, i}^{t} - \bLa_{8, i}^{t} - \bLa_{10, i}^{t} + \rho\bu_i^{t}\mathbf{1}^{t \top} + \rho\bW_i^{t} + \rho\bH_i^{t} + \rho\bZ_i^{t}) (\rho\mathbf{1}\mathbf{1}^{\top} + 3\rho\bI)^{-1}$ \label{step:S}
	    \State $\beta^{t+1}_i = \Pi\left(\frac{1}{\rho}(-\bla_{6, i}^{t} - \rho\bH_i^{t \top}\bdelta + \rho\bb_i\odot(\bB_i\boldsymbol{\eta}))\right)_{\bbR_{+}}$ \label{step:beta}
	    \EndFor
     
	    \State $\tilde{\bc}^{t+1} = \Pi(\frac{1}{\rho} (\tilde{\bI}^{\top}\tilde{\bI} + \tilde{\mathbf{1}} \tilde{\mathbf{1}}^{\top})^{-1} (\tilde{\bI}^{\top}\bla_7^{t} - \lambda_9^{t} \tilde{\mathbf{1}} - \rho \tilde{\bI}^{\top} (\balpha \odot \bgamma) $
         \Statex $+ \rho \tilde{\bI}^{\top} (\balpha \odot (\bDelta^{\top} \bu^{t})) - \rho \tilde{\bbeta} \tilde{\mathbf{1}} + \rho \Omega \tilde{\mathbf{1}} )_{\bbR_{+}}$ \label{step:c}
     
	    \State $\bu^{t+1} = \frac{1}{\rho}(\bI + \tilde{\bR}^{\top}\tilde{\bR} + \tilde{\bD}^{\top}\tilde{\bD} + (\Delta\tilde{\alpha})(\Delta\tilde{\alpha})^{\top})^{-1} (\bla_{1}^{t} + \tilde{\bR}^{\top}\bla_2^{t} - \tilde{\bD}^{\top}\bla_3^{t} - (\bDelta\tilde{\balpha})\bla_7^{t} + \rho\tilde{\bu}^{t+1} - \rho\tilde{\bR}^{\top}\omega^{t+1} + \rho\tilde{\bD}^{\top}\bq + \rho(\bDelta\tilde{\balpha})(\balpha \odot \bgamma) +\rho(\bDelta\tilde{\balpha})(\tilde{\mathbf{I}}\tilde{\bc}^{t+1}))$  \label{step:u}
	    
	    \For {$i = 1, \ldots, n$}
        \State $\bW_i^{t+1} = \frac{1}{\rho}(\mathbf{1}\mathbf{1}^{\top} + \bI)^{-1} (\rho\mathbf{1}\mathbf{1}^{\top} + \rho\bS_i^{t+1} - \mathbf{1}\bla_{4, i}^{t \top} + \bLa_{5, i}^{t} ) $ \label{step:W}
        
        \State $\bH_i^{t+1} = \frac{1}{\rho}(\bdelta \bdelta^{\top} + \bI)^{-1} (-\bdelta \bla_{6, i}^{t \top} + \bLa_{8, i}^{t} - \rho \bdelta \beta_{i}^{t+1 \top} + \rho \bdelta(\bb_{i} \odot \bB_{i}\boldsymbol{\eta})^{\top} + \rho \bS_{i}^{t+1})$
        \label{step:H}
        
	   \State $\bZ_i^{t+1} = \mathbbm{1}({\rho>\tilde{\lambda}}) \Pi {\left(\left(\frac{1}{\rho-\tilde{\lambda}}\right) (\rho\bS_{i}^{t+1} + \bLa_{10}^{t} - \frac{\tilde{\lambda}}{2})\right)}_{[0, 1]} + \mathbbm{1}({\rho<\tilde{\lambda}}) \Pi {\left(\left(\frac{1}{\rho-\tilde{\lambda}}\right)(\rho\bS_i^{t+1} + \bLa_{10}^{t} - \frac{\tilde{\lambda}}{2})\right)}_{\{0, 1\}}$ \label{step:Z}
	    \EndFor
	    
	    \For {$i = 1, \ldots, n$}
	    \State $\bla^{t+1}_{1, i} = \bla^{t}_{1, i} + \rho (\bS_i^{t+1}\mathbf{1} - \bu_i^{t+1})$  \label{step:lambda1}
	    
	    \State $\bla^{t+1}_{3, i} = \bla^{t}_{3, i} + \rho (\bD \bu_i^{t+1} - \bq_i)$ \label{step:lambda3}
	    	    
	    \State $\bla^{t+1}_{4, i} = \bla^{t}_{4, i} + \rho (\bW_i^{t+1 \top} \mathbf{1} -\mathbf{1})$ \label{step:lambda4}
	    	    
	    \State $\bLa^{t+1}_{5, i} = \bLa^{t}_{5, i} + \rho (\bS_i^{t+1}  - \bW_i^{t+1})$ \label{step:lambda5}
	    	    
	    \State $\bla^{t+1}_{6, i} = \bla^{t}_{6, i} + \rho (\bH_i^{t+1 \top} \bdelta + \beta_i^{t+1} - \bb_i \odot \bB_i \boldsymbol{\eta})$ \label{step:lambda6}
	    	    
	    \State $\bLa^{t+1}_{8, i} = \bLa^{t}_{8, i} + \rho (\bS_i^{t+1} - \bW_i^{t+1})$ \label{step:lambda8}
	    	    
	    \State $\bLa^{t+1}_{10, i} = \bLa^{t}_{10, i} + \rho (\bS_i^{t+1} - \bZ_i^{t+1})$ \label{step:lambda10}
	    \EndFor
	    	    
	    \State $\bla^{t+1}_{2} = \bla^{t}_{2} + \rho (\bomega^{t+1} - \bR(\sum_{i=1}^n \bu_i^{t+1}))$ \label{step:lambda2}
	    
	    \State $\bla^{t+1}_{7} = \bla^{t}_{7} + \rho (\balpha \odot (\bDelta^{\top}\bu^{t+1} - \bdelta) - \tilde{\bI}\tilde{\bc}^{t+1})$ \label{step:lambda7}
	    
	    \State $\lambda^{t+1}_{9} = \lambda^{t}_{9} + \rho (\tilde{\bc}^{t+1 \top}\tilde{\mathbf{1}} + \tilde{\bbeta}^{t+1} - \Omega)$ \label{step:lambda9}
	    
		\EndFor
		\State \textbf{Return:} $\bS_i^{\tilde{T}}, \forall i= 1, \dots, n $
	\end{algorithmic}
\normalsize
\end{algorithm}
\normalsize
Algorithm~\ref{alg:ADMM-ForOurProblem} solves the relaxed version of problem~\eqref{eq:optModel1}. In this algorithm, we use the projection operator $\Pi(\cdot)_{[0,1]}$  that projects elements of a matrix to the interval $[0,1]$. $\Pi(\cdot)_{\bbR_{+}}$ is also a projection operator but projects elements of a matrix to $\bbR_{+}$. Notice that in Algorithm~\ref{alg:ADMM-ForOurProblem}, the computation load of steps~\ref{step:S}, \ref{step:W}, \ref{step:H}, and \ref{step:Z} is extensive because matrices $\bS, \bW, \bH$ and $\bZ$ have large sizes. However, each column in these matrices corresponds to one driver and these steps are not coupled so we can perform the computation of each column in parallel by leveraging parallel computation. 
The notations used in Algorithm~\ref{alg:ADMM-ForOurProblem} are defined below.
\vspace{-0.6cm}
\begin{figure}[H]
  \begin{minipage}[b]{0.13\textwidth}
  \[
              \bgamma =
              \begin{bmatrix}
                \gamma_1 \\
                \vdots\\
                \gamma_n 
              \end{bmatrix},
            \]
  \end{minipage}
  \centering
  \begin{minipage}[b]{0.13\textwidth}
  \[
              \bq =
              \begin{bmatrix}
                \bq_1 \\
                \vdots\\
                \bq_n 
              \end{bmatrix},
            \]
  \end{minipage}
  \centering
  \begin{minipage}[b]{0.24\textwidth}
  \[
              \bla_i =
              \begin{bmatrix}
                \bla_{i, 1} \\
                \vdots\\
                \bla_{i, n} 
              \end{bmatrix}, i=1,3,
            \]
  \end{minipage}
  \centering
  \begin{minipage}[b]{0.12\textwidth}
  \[
              \balpha =
              \begin{bmatrix}
                \balpha_1\\
                \vdots\\
                \balpha_n
              \end{bmatrix},
            \]
  \end{minipage}
  \centering
  \begin{minipage}[b]{0.14\textwidth}
  \[
              \tilde{\bu}^{t} =
              \begin{bmatrix}
                \bS_1^t \mathbf{1} \\
                \vdots\\
                \bS_n^t \mathbf{1}
              \end{bmatrix},
            \]
  \end{minipage}
  \centering
  \begin{minipage}[b]{0.1\textwidth}
  \[
          \tilde{\mathbf{1}} =
          \begin{bmatrix}
            \mathbf{1}\\
            \mathbf{0}
          \end{bmatrix},
        \]
  \end{minipage}
  \centering
  \begin{minipage}[b]{0.1\textwidth}
  \[
          \tilde{\bc} =
          \begin{bmatrix}
            \bc\\
            \bmu
          \end{bmatrix},
        \]
  \end{minipage}
  \normalsize   
\end{figure}  
\vspace{-0.9cm}
\normalsize   
\begin{figure}[H]
  \centering
  \begin{minipage}[b]{0.2\textwidth}
  \[
              \tilde{\balpha} =
              \begin{bmatrix}
                \alpha_1 & & \\
                & \ddots & \\
                & & \alpha_n
              \end{bmatrix},
            \]
  \end{minipage}
  \centering
  \begin{minipage}[b]{0.2\textwidth}
  \[
              \tilde{\bD} =
              \begin{bmatrix}
                \bD & & \\
                & \ddots & \\
                & & \bD
              \end{bmatrix},
            \]
  \end{minipage}
  \centering
  \begin{minipage}[b]{0.2\textwidth}
  \[
              \bDelta =
              \begin{bmatrix}
                \bdelta & & \\
                & \ddots & \\
                & & \bdelta
              \end{bmatrix},
            \]
  \end{minipage}
  \centering
  \begin{minipage}[b]{0.2\textwidth}
  \[
          \tilde{\bR} =
          \begin{bmatrix}
            \bR & \dots & \bR 
          \end{bmatrix},
        \]
  \end{minipage}
  \centering
  \begin{minipage}[b]{0.15\textwidth}
  \[
          \tilde{\bI} =
          \begin{bmatrix}
            \bI & -\bI  
          \end{bmatrix}.
        \]
  \end{minipage}
    \centering
\end{figure}
\normalsize

\vspace{0.2cm}

\section{Reformulated Optimization Model for the ADMM Algorithm} \label{appndx:optModelReformulatedSection}
To solve the relaxed version of problem~\eqref{eq:optModel1} efficiently, we present a distributed algorithm based on this reformulation

\begingroup
\allowdisplaybreaks
\begin{equation}   
\label{eq:optModelReformulated}
    \begin{split}
        \min_{\substack{\mathbf{S}, \mathbf{H}, \mathbf{W}, \mathbf{Z}, \mathbf{u}, \boldsymbol{\beta},\\
        \omega, \mu_i, \tilde{\beta}, \mathbf{c}}} \quad & \sum_{\ell=1}^{|\cE|} \sum_{t=1}^{|\mathbf{T}|} \hat{v}_{\ell, t} \theta_{\ell, t}(\hat{v}{\ell, t}) \\ &- \frac{\tilde{\lambda}}{2} \sum_{r=1}^{\cR} \sum_{t=1}^{|\bT|} \sum_{i=1}^{n}(\bZ_i)_{r, t}((\bZ_i)_{r, t} - 1)\\
        \stt \quad & \bS_i \mathbf{1} = \bu_i,\quad \forall i=1, 2, \dots, n\\
        & \bomega = \tilde{\bR}\bu \\
        & \tilde{\bD} \bu = \bq,\quad \forall i=1, 2, \dots, n \\
        & \bW_i^{\top} \mathbf{1} = \mathbf{1},\quad \forall i=1, 2, \dots, n \\
        & \bS_i = \bW_i,\quad \forall i=1, 2, \dots, n \\
        & \bH_i^{\top} \bdelta + \bbeta_i = \bb_i \odot \bB_i \boldsymbol{\eta},\quad \forall i=1, 2, \dots, n\\
        & \bS_i = \bH_i,\quad \forall i=1, 2, \dots, n \\
        & \bbeta_i \geq 0,\quad \forall i=1, 2, \dots, n\\ 
        & \bZ_i \in \lbrack0, 1\rbrack^{(|\cR|\cdot|\bT|)\times(|\cN_{i}|)}, \quad \forall i=1, 2, \dots, n \\
        & \tilde{\bI} \tilde{\bc} = \balpha \odot (\bDelta \bu - \bgamma) \\
        & \tilde{\bc} \geq 0, \quad \tilde{\bc}^{\top} \tilde{\mathbf{1}} + \tilde{\beta} = \Omega, \quad \tilde{\beta} \geq 0 \\
        & \bS_i = \bZ_i,\quad \forall i=1, 2, \dots, n,
    \end{split}
    \normalsize
\end{equation}
\endgroup
\normalsize
where $\bS = \{\bS_i\}_{i=1}^n$, $\bH = \{\bH_i\}_{i=1}^n$, $\bW = \{\bW_i\}_{i=1}^n$, $\bZ = \{\bZ_i\}_{i=1}^n$, $\bu = \{\bu_i\}_{i=1}^n$, and $\boldsymbol{\beta} = \{\boldsymbol{\beta}_i\}_{i=1}^n$.

\vspace{0.2cm}

\section{Details of Alternating Direction Method of Multipliers (ADMM)}
\label{appdx:ADMM}
Before explaining the steps of our proposed algorithm, let us first explain the Alternating Direction Method of Multipliers (ADMM), which is a main building block of our framework.
\subsection{Review of ADMM} \label{appdx:ADMMReview}
ADMM developed in~\cite{gabay1976dual} and \cite{glowinski1975approximation} aims at solving linearly constrained optimization problems of the form
\begin{align}
    \min_{w,z} \;\;  h(w) + g(z)\quad \quad 
    \stt \quad  Aw + Bz = c, \nonumber
\normalsize
\end{align}
\normalsize
where $w\in \mathbb{R}^{d_1}, z\in \mathbb{R}^{d_2}$, $c \in \mathbb{R}^{k}$, $A \in \mathbb{R}^{k \times d_1}$, and $B \in \mathbb{R}^{k \times d_2}$. By forming the augmented Lagrangian function 
\begin{equation*}
\mathcal{L}(w,z,\lambda) \triangleq h(w) + g(z) + \langle \lambda , Aw + Bz - c\rangle  + \frac{\rho}{2} \|Aw + Bz - c\|_2^2,
\normalsize
\end{equation*}
\normalsize
each iteration of ADMM applies  alternating minimization to the primal variables and  gradient ascent  to the dual variables:
\begin{align}
    & \textrm{Primal Update:} &&w^{r+1}  = \arg\min_{w}\; \mathcal{L}(w,z^r,\lambda^r),\\
    & &&z^{r+1}  = \arg\min_{z}\; \mathcal{L}(w^{r+1},z,\lambda^r)\nonumber\\
    & \textrm{Dual Update:} &&\lambda^{r+1}  = \lambda^r + \rho \left( Aw^{r+1} + Bz^{r+1} - c\right)  \nonumber
\normalsize
\end{align}
\normalsize
This algorithm is extensively explored in the optimization literature (see \cite{boyd2011distributed} for a monograph on the use of this algorithm in convex distributed optimization and \cite{hong2016convergence} for its use in non-convex continuous optimization).

\subsection{ADMM for Solving  the Relaxed Version of Problem~\eqref{eq:optModel1}}
\label{appdx:ADMMSteps}
Let
{\begingroup
\allowdisplaybreaks
\begin{align}
\allowdisplaybreaks
    & \mathcal{L}(\{\bS_{i}\}_{i=1}^{n}, \{\bH_{i}\}_{i=1}^{n}, \{\bW_{i}\}_{i=1}^{n}, \{\bZ_{i}\}_{i=1}^{n}, \bomega, \tilde{\beta}, \{\bbeta_{i}\}_{i=1}^{n},  \tilde{\bc},\{\bu_{i}\}_{i=1}^{n})  \nonumber\\
    \triangleq 
    & \; F(\bomega) - \frac{\tilde{\lambda}}{2} \sum_{i=1}^{n} \sum_{j=1}^{|\cN_i|} \sum_{\ell=1}^{|\cE|} \sum_{t=1}^{|\mathbf{T}|} (\bZ_j)_{i, (r, t)}((\bZ_j)_{i, (r, t)} - 1) \nonumber \\
    & + \sum_{i=1}^{n} \mathbb{I}_{\mathbb{R}_{+}^{|\cN_i|}}(\bbeta_i) + \mathbb{I}_{\mathbb{R}_{+}^{2n}}(\tilde{\bc})  + \mathbb{I}_{\mathbb{R}_{+}}(\tilde{\beta}) + \mathbb{I}_{[0, 1]}((\bZ_j)_{i, (r, t)}) \nonumber \\
    & + \sum_{i=1}^{n} \langle \bla_{1, i}, \bS_i \mathbf{1} - \bu_i \rangle + \langle\bla_2, \bomega - \tilde{\bR}\bu\rangle  \nonumber\\
    & + \sum_{i=1}^{n}\langle\bla_{3, i}, \bD\bu_i-\bq_i\rangle  + \sum_{i=1}^{n}\langle\bla_{4, i}, \bW_i\mathbf{1}-\mathbf{1}\rangle \nonumber\\
    & + \sum_{i=1}^{n}\langle\bLa_{5, i},\bS_i-\bW_i\rangle + \sum_{i=1}^{n} \langle \bla_{6, i}, \bH_i^{\top} \delta_p + \bbeta_i - \bb_i \cdot (\bB_i \boldsymbol{\eta}) \rangle \nonumber\\
    & +  \langle \bla_7, (\Delta \tilde{\balpha})^{\top} \bu - \tilde{\balpha} \bdelta - \tilde{\bI} \tilde{\bc} \rangle + \sum_{i=1}^{n}\langle\bLa_{8, i}, \bS_i-\bH_i\rangle \nonumber \\
    &   +  \langle\bla_9,\tilde{\bc}^{\top} \tilde{\mathbf{1}} + \tilde{\beta} - \Omega\rangle +\langle \bLa_{10,i},\bS_i-\bZ_i\rangle  \nonumber \\
    & + \frac{\rho}{2}\sum_{i=1}^{n}||\bS_i \mathbf{1} - \bu_i||^{2} + \frac{\rho}{2} \sum_{i=1}^{n} ||\bS_i \mathbf{1} - \bu_i||^{2} + \frac{\rho}{2} ||\bomega - \tilde{\bR} \bu||^{2} \nonumber \\
    & + \frac{\rho}{2}\sum_{i=1}^{n}||\bD\bu_i-\bq_i||^{2} + \frac{\rho}{2}\sum_{i=1}^{n}||\bW_i\mathbf{1}-\mathbf{1}||^{2} \nonumber\\
    & +  \frac{\rho}{2}\sum_{i=1}^{n}||\bS_i-\bW_i||^{2} + \frac{\rho}{2}\sum_{i=1}^{n}||\bH_i^{\top}\bdelta+\bbeta_i - \bb_i \cdot (\bB_i\boldsymbol{\eta})||^{2} \nonumber \\
    &   + \frac{\rho}{2}||(\Delta\tilde{\balpha})^{\top}\bu-\tilde{\balpha}\bdelta-\tilde{\bI}\tilde{\bc}||^{2} + \frac{\rho}{2}\sum_{i=1}^{n}||\bS_i-\bH_i||^{2}  \nonumber \\
    &  + \frac{\rho}{2}||\tilde{\bc}^{\top}\mathbf{1}+\tilde{\beta}-\Omega||^{2} + \frac{\rho}{2}\sum_{i=1}^{n}||\bS_i-\bZ_i||^{2} \nonumber \\
\normalsize
\end{align}
\normalsize
\endgroup}
\normalsize
\noindent{\normalsize be the augmented Lagrangian function of the relaxed version of problem~\eqref{eq:optModel1} with the set of Lagrange multipliers $\{\{\bla_1\}_{i=1}^{n}, \bla_2, \dots, \{\bLa_{10}\}_{i=1}^{n}\}$ and $\rho > 0$ be the primal penalty parameter. Then, ADMM solves the relaxed version of problem~\eqref{eq:optModel1} by the following iterative scheme}
\begingroup
\allowdisplaybreaks
\begin{align*}
\bomega_{(\ell, \hat{t})}^{t+1} =  
& \argmin\limits_{\bomega_{(\ell, \hat{t})}} \quad  \bomega_{(\ell, \hat{t})} \theta_{\ell, t}(\bomega_{(\ell, \hat{t})}) \\
& + (\bla_2^{t})_{\ell, \hat{t}}(\bomega_{\ell, \hat{t}} - \br_{\ell, \hat{t}} (\sum_{i=1}^{n}\bu_i)) + \frac{\rho}{2} (\bomega_{\ell, \hat{t}} - \br_{\ell, \hat{t}} (\sum_{i=1}^{n}\bu_i^{t})) \\ 
\bS_i^{t+1} = 
& \argmin\limits_{\bS_i} \quad   \langle \bla_{1, i}^{t}, \bS_i \mathbf{1} - \bu_i^{t} \rangle + \langle\bLa_{5, i}^{t},\bS_i-\bW_i^{t}\rangle \\
& + \langle\bLa_{8, i}^{t}, \bS_i-\bH_i^{t}\rangle +\langle \bLa_{10,i}^{t},\bS_i-\bZ_i^{t}\rangle \\
& + \frac{\rho}{2}||\bS_i \mathbf{1} - \bu_i^{t}||^{2} + \frac{\rho}{2} \sum_{i=1}^{n} ||\bS_i \mathbf{1} - \bu_i^{t}||^{2} \\
& + \frac{\rho}{2}||\bS_i-\bW_i^{t}||^{2}  + \frac{\rho}{2}||\bS_i-\bH_i^{t}||^{2} + \frac{\rho}{2}||\bS_i-\bZ_i^{t}||^{2} \\
& ,\forall i=1, 2, \dots, n \\
\bbeta_i^{t+1} = & \argmin\limits_{\bbeta_i} \quad  \mathbb{I}_{\mathbb{R}_{+}^{|\cN_i|}}(\bbeta_i) \\
& + \langle \bla_{6, i}^{t}, \bH_i^{t \top} \delta_p + \bbeta_i^{t} - \bb_i \cdot (\bB_i \boldsymbol{\eta}) \rangle \\
& + \frac{\rho}{2}||\bH_i^{t \top} \bdelta + \bbeta_i^{t} - \bb_i\cdot\bB_i\boldsymbol{\eta}||^{2}, \quad \forall i=1, 2, \dots, n \\
\tilde{\bc}^{t+1} = & \argmin\limits_{\tilde{\bc}} \quad \mathbb{I}_{\mathbb{R}_{+}^{2n}}(\tilde{\bc}) + \langle \bla_7^{t}, (\Delta \tilde{\balpha})^{\top} \bu^{t} - \tilde{\balpha} \bdelta - \tilde{\bI} \tilde{\bc} \rangle \\
& + \langle\bla_9^{t},\tilde{\bc}^{\top} \tilde{\mathbf{1}} + \tilde{\beta} - \Omega\rangle + \frac{\rho}{2}||(\Delta\tilde{\balpha})^{\top}\bu^{t} - \tilde{\balpha}\bdelta - \tilde{\bI}\tilde{\bc}||^{2} \\
&  + \frac{\rho}{2}||\tilde{\bc}^{\top} \mathbf{1} + \tilde{\beta}^{t} - \Omega||^{2} \\
\bu^{t+1} = & \argmin\limits_{\bu} \quad \langle \bla_1^{t}, \tilde{\bu}^{t+1} - \bu \rangle + \langle\bla_2^{t}, \bomega - \tilde{\bR}\bu\rangle  \\
& + \langle\bla_{3}^{t}, \tilde{\bD}\bu -\bq\rangle + \langle \bla_7^{t}, (\Delta \tilde{\balpha})^{\top} \bu - \tilde{\balpha} \bdelta - \tilde{\bI} \tilde{\bc}^{t+1} \rangle  \\
&  + \frac{\rho}{2} ||\tilde{\bu}^{t+1} - \bu||^{2} + \frac{\rho}{2} ||\bomega - \tilde{\bR} \bu||^{2} + \frac{\rho}{2}||\tilde{\bD}\bu- \bq||^{2} \\
&  + \frac{\rho}{2}||(\Delta\tilde{\balpha})^{\top}\bu-\tilde{\balpha}\bdelta-\tilde{\bI}\tilde{\bc}^{t+1}||^{2} \\
\bW_i^{t+1} = & \argmin\limits_{\bW_i} \quad \langle\bla_{4, i}^{t}, \bW_i\mathbf{1} - \mathbf{1}\rangle + \langle\bLa_{5, i}^{t}, \bS_i^{t+1}-\bW_i\rangle \\
& + \frac{\rho}{2} ||\bW_i\mathbf{1} - \mathbf{1}||^{2} + \frac{\rho}{2}||\bS_i^{t+1}-\bW_i||^{2}\\
& , \quad \forall i=1, 2, \dots, n\\
\bH_i^{t+1} = & \argmin\limits_{\bH_i} \quad \langle \bla_{6, i}^{t}, \bH_i^{\top} \delta_p + \bbeta_i^{t+1} - \bb_i \cdot (\bB_i \boldsymbol{\eta}) \rangle \\
& + \langle\bLa_{8, i}^{t}, \bS_i^{t+1}-\bH_i\rangle \\
& + \frac{\rho}{2} ||\bH_i^{\top}\bdelta + \bbeta_i^{t} - \bb_i \cdot (\bB_i\boldsymbol{\eta})||^{2} \\
& + \frac{\rho}{2} ||\bS_i^{t+1}-\bH_i||^{2}, \quad \forall i=1, 2, \dots, n \\
\bZ_i^{t+1} = & \argmin\limits_{\bZ_i} \quad
\mathbbm{1}({\rho>\tilde{\lambda}}) \mathbb{I}_{[0, 1]^{(|\cR|\cdot|\bT|) \times |\cN_i|}}(\bZ_i) \\
& + \mathbbm{1}({\rho<\tilde{\lambda}})\mathbb{I}_{\{0, 1\}^{(|\cR|\cdot|\bT|) \times |\cN_i|}}(\bZ_i) \\
& - \frac{\tilde{\lambda}}{2} \sum_{j=1}^{|\cN_i|} \sum_{\ell=1}^{|\cE|} \sum_{t=1}^{|\mathbf{T}|} (\bZ_j)_{i, (r, t)}((\bZ_j)_{i, (r, t)} - 1) \\
& + \langle \bLa_{10,i}^{t},\bS_i^{t+1}-\bZ_i\rangle + \frac{\rho}{2}\sum_{i=1}^{n}||\bS_i^{t+1}-\bZ_i||^{2}\\
& , \quad \forall i=1, 2, \dots, n\\
\tilde{\beta}^{t+1} = & \argmin\limits_{\tilde{\beta}} \quad \mathbb{I}_{\mathbb{R}_{+}}(\tilde{\beta}) \\
& +  \langle\bla_9^{t},\tilde{\bc}^{t+1 \top} \tilde{\mathbf{1}} + \tilde{\beta} - \Omega\rangle \\
& + \frac{\rho}{2}||\tilde{\bc}^{t+1 \top}\mathbf{1} + \tilde{\beta} - \Omega||^{2}\\
\bla_{1,i}^{t+1} = & \bla_{1,i}^{t} + \rho(\bS_i^{t+1}\mathbf{1} - \bu_i^{t+1}), \forall i=1, 2, \dots, n \\
\bla_{2}^{t+1} = & \bla_{2}^{t} + \rho(\bomega^{t+1} - \bR(\sum_{i=1}^{n}\bu_i^{t+1})) \\
\bla_{3,i}^{t+1} = & \bla_{3,i}^{t} + \rho(\bD\bu_i^{t+1} - \bq_i), \forall i=1, 2, \dots, n \\
\bla_{4,i}^{t+1} = & \bla_{4,i}^{t} + \rho(\bW_i^{t+1 \top}\mathbf{1} - \mathbf{1}), \forall i=1, 2, \dots, n \\
\bLa_{5,i}^{t+1} = & \bLa_{5,i}^{t} + \rho(\bS_i^{t+1} - \bW_i^{t+1}), \forall i=1, 2, \dots, n \\
\bla_{6,i}^{t+1} = & \bla_{6,i}^{t} + \rho(\bH_i^{t+1 \top}\bdelta + \bbeta_i^{t+1} - \bb_i \odot \bB_i \boldsymbol{\eta}),\\ \forall i=1, 2, \dots, n \\
\bla_{7}^{t+1} = & \bla_{7}^{t} + \rho(\balpha \odot (\Delta^{\top}\bu^{t+1} - \bdelta) - \tilde{\bI}\tilde{\bc}^{t+1}) \\
\bLa_{8,i}^{t+1} = & \bLa_{8,i}^{t} + \rho(\bS_i^{t+1} - \bH_i^{t+1}), \forall i=1, 2, \dots, n \\
\bla_{9}^{t+1} = & \bla_{9}^{t} + \rho(\tilde{\bc}^{t+1 \top}\mathbf{1} + \tilde{\beta}^{t+1} - \Omega) \\
\bLa_{10,i}^{t+1} = & \bLa_{10,i}^{t} + \rho(\bS_i^{t+1} - \bZ_i^{t+1}), \forall i=1, 2, \dots, n\\
\end{align*}
\endgroup
\normalsize
The primal update rules can be simplified as
\begingroup
\allowdisplaybreaks
\begin{align*}
\omega^{t+1}_{\ell, \hat{t}} = 
& \argmin\limits_{\omega_{\ell, \hat{t}}}   \omega_{\ell, \hat{t}} \theta_{\ell, t}(\omega_{\ell, \hat{t}}) \\
& + \bla^{t}_{2, (\ell, \hat{t})}(\omega_{\ell, \hat{t}} - \br_{\ell, \hat{t}} \left(\sum_{i=1}^{n}\bu_i^{t}\right)) \\
& + \frac{\rho}{2} (\omega_{\ell, \hat{t}} - \bR_{\ell, \hat{t}} \left(\sum_{i=1}^{n}\bu_i^{t}\right))^{2}\\
&, \quad \forall \ell=1, 2, \dots, |\cE|, \forall \hat{t}=1, 2, \dots, \tilde{T} \\
\bS^{t+1}_i = & (-\bla_{1,i}^{t}\mathbf{1}^{\top} - \bLa_{5, i}^{t} - \bLa_{8, i}^{t} - \bLa_{10, i}^{t} \\
& + \rho\bu_i\mathbf{1}^{t \top} + \rho\bW_i^t + \rho\bH_i^{t} + \rho\bZ_i^{t}) (\rho\mathbf{1}\mathbf{1}^{\top} + 3\rho\bI)^{-1}\\
&, \forall i =1, 2, \dots, n \\
\beta^{t+1}_i = & \Pi\left(\frac{1}{\rho}(-\bla_{6, i}^t - \rho\bH_i^{t \top}\bdelta + \rho\bb_i\odot(\bB_i\boldsymbol{\eta}))\right)_{\bbR_{+}}\\
&, \forall i =1, 2, \dots, n \\
\tilde{\bc}^{t+1} = & \Pi(\frac{1}{\rho} (\tilde{\bI}^{\top}\tilde{\bI} + \tilde{\mathbf{1}} \tilde{\mathbf{1}}^{\top})^{-1} (\tilde{\bI}^{\top}\bla_7^t - \lambda_9^t \tilde{\mathbf{1}} - \rho \tilde{\bI}^{\top} (\balpha \odot \bgamma)\\
&  + \rho \tilde{\bI}^{\top} (\balpha \odot (\bDelta^{\top} \bu^t)) - \rho \tilde{\bbeta} \tilde{\mathbf{1}} + \rho \Omega \tilde{\mathbf{1}})_{\bbR_{+}} \\
\bu^{t+1} = & \frac{1}{\rho}(\bI + \tilde{\bR}^{\top}\tilde{\bR} + \tilde{\bD}^{\top}\tilde{\bD} + (\Delta\tilde{\alpha})(\Delta\tilde{\alpha})^{\top})^{-1}\\
&(\bla_{1}^{t} + \tilde{\bR}^{\top}\bla_2^{t} - \tilde{\bD}^{\top}\bla_3^{t} - (\bDelta\tilde{\balpha})\bla_7^{t} + \rho\tilde{\bu}^{t} \\
& - \rho\tilde{\bR}^{\top}\omega^{t+1} + \rho\tilde{\bD}^{\top}\bq + \rho(\bDelta\tilde{\balpha})(\balpha \odot \bgamma)\\
& +\rho(\bDelta\tilde{\balpha})(\tilde{\mathbf{I}}\tilde{\bc}^{t+1})) \\
\bZ_i^{t+1} = & \mathbbm{1}({\rho>\tilde{\lambda}}) \Pi {\left(\left(\frac{1}{\rho-\tilde{\lambda}}\right) (\rho\bS_{i}^{t+1} + \bLa_{10}^{t} - \frac{\tilde{\lambda}}{2})\right)}_{[0, 1]} \\
& + \mathbbm{1}({\rho<\tilde{\lambda}}) \Pi {\left(\left(\frac{1}{\rho-\tilde{\lambda}}\right)(\rho\bS_i^{t+1} + \bLa_{10}^{t} - \frac{\tilde{\lambda}}{2})\right)}_{\{0, 1\}}\\
&, \forall i =1, 2, \dots, n \\
\bW_i^{t+1} = & \frac{1}{\rho}(\mathbf{1}\mathbf{1}^{\top} + \bI)^{-1} (-\mathbf{1}\bla_{4, i}^{t \top} + \bLa_{5, i}^t + \rho\mathbf{1}\mathbf{1}^{\top} + \rho\bS_i^{t+1})\\
&, \forall i =1, 2, \dots, n \\
\bH_i^{t+1} = & \frac{1}{\rho} (\bdelta\bdelta^{\top} + \bI)^{\top}(-\bdelta\bla_{6, i}^{t \top} + \bLa_{8, i}^t - \rho \bdelta \bbeta_{i}^{\top} \\
& + \rho \bdelta(\bb_i \cdot \bB_i \boldsymbol{\eta})^{\top} + \rho\bS_i^{t+1}), \forall i =1, 2, \dots, n \\  
\end{align*}
\normalsize
\endgroup
\normalsize

\end{document}